# A BOXING INEQUALITY FOR THE FRACTIONAL PERIMETER

AUGUSTO C. PONCE AND DANIEL SPECTOR

ABSTRACT. We prove the Boxing inequality:
$$\mathcal{H}_\infty^{d-\alpha}(U) \leq C\alpha(1-\alpha) \int_U \int_{\mathbb{R}^d \setminus U} \frac{\mathrm{d}y\,\mathrm{d}z}{|y-z|^{\alpha+d}},$$
for every $\alpha \in (0,1)$ and every bounded open subset $U \subset \mathbb{R}^d$, where $\mathcal{H}_\infty^{d-\alpha}(U)$ is the Hausdorff content of $U$ of dimension $d-\alpha$ and the constant $C > 0$ depends only on $d$. We then show how this estimate implies a trace inequality in the fractional Sobolev space $W^{\alpha,1}(\mathbb{R}^d)$ that includes Sobolev's $L^{\frac{d}{d-\alpha}}$ embedding, its Lorentz-space improvement, and Hardy's inequality. All these estimates are thus obtained with the appropriate asymptotics as $\alpha$ tends to 0 and 1, recovering in particular the classical inequalities of first order. Their counterparts in the full range $\alpha \in (0,d)$ are also investigated.

## 1. INTRODUCTION AND MAIN RESULTS

Let $d \in \mathbb{N}$ and write $\mathbb{R}^d$ to denote the Euclidean space of $d$ dimensions. A geometric formulation of the classical Boxing inequality of W. Gustin [26] can be stated as

**Theorem 1.1.** *There exists a constant $C = C(d) > 0$ such that for every bounded open set $U \subset \mathbb{R}^d$ with smooth boundary one can find a covering*
$$U \subset \bigcup_{i=0}^\infty B_{r_i}(x_i)$$
*by open balls of radii $r_i$ for which*
$$\sum_{i=0}^\infty r_i^{d-1} \leq C \operatorname{Per}(U).$$

Here we utilize $\operatorname{Per}(U)$ to denote the perimeter of $U$, i.e. integration of the $(d-1)$-dimensional measure over the topological boundary $\partial U$, since this inequality has been shown to hold in the more general class of sets of finite perimeter in the sense of De Giorgi; see e.g. [20, Corollary 4.5.4].

---

One can think of such an estimate as a $(d-1)$-dimensional analogue of the trivial fact that for every bounded open set $U \subset \mathbb{R}^d$ one can find a covering by open balls of radii $r_i$ for which

$$\sum_{i=0}^{\infty} \omega_d r_i^d \leq 2|U|,$$

where $\omega_d$ denotes the volume of the unit ball $B_1 \subset \mathbb{R}^d$ and $|U|$ is the Lebesgue measure of the set $U$. The constant $2$ is only for convenience of display, and can be replaced by any $C > 1$, as the inequality merely follows from the definition of the Lebesgue measure.

The principal new result of this paper is the following theorem that interpolates these two estimates,

**Theorem 1.2.** *There exists a constant $C = C(d) > 0$ such that for every bounded open set $U \subset \mathbb{R}^d$ one can find a covering*

$$U \subset \bigcup_{i=0}^{\infty} B_{r_i}(x_i)$$

*by open balls of radii $r_i$ for which*

$$\sum_{i=0}^{\infty} r_i^{d-\alpha} \leq C\alpha(1-\alpha) P_\alpha(U),$$

*for every $\alpha \in (0,1)$.*

Here $P_\alpha$ is defined for bounded open sets and more generally for any Lebesgue measurable set $A \subset \mathbb{R}^d$ by

(1.1) $$P_\alpha(A) := 2 \int_A \int_{\mathbb{R}^d \setminus A} \frac{\mathrm{d}y \, \mathrm{d}z}{|y-z|^{\alpha+d}}.$$

It has been called the fractional perimeter [16] or non-local $\alpha$-perimeter [14] and gives one notion of an intermediate object between the classical perimeter and the Lebesgue measure. One observes that $P_\alpha$ enjoys an isoperimetric inequality (see [5, 25] and also [23])

$$\frac{P_\alpha(B_1)}{|B_1|^{\frac{d-\alpha}{d}}} \leq \frac{P_\alpha(A)}{|A|^{\frac{d-\alpha}{d}}},$$

as well as a coarea formula for functions in the fractional Sobolev space $W^{\alpha,1}(\mathbb{R}^d)$:

(1.2) $$[u]_{W^{\alpha,1}(\mathbb{R}^d)} := \int_{\mathbb{R}^d} \int_{\mathbb{R}^d} \frac{|u(y) - u(z)|}{|y-z|^{\alpha+d}} \, \mathrm{d}y \, \mathrm{d}z = \int_{-\infty}^{\infty} P_\alpha(\{u > t\}) \, \mathrm{d}t;$$

see [43] and Lemma 4.3 below. Moreover, one has the asymptotics

(1.3) $$\lim_{\alpha \to 0} \alpha P_\alpha(A) = C'|A| \quad \text{and} \quad \lim_{\alpha \to 1} (1-\alpha) P_\alpha(A) = C'' \operatorname{Per}(A),$$

that allows recovery of the endpoints; see [18, 32, 41].

Our proof of Theorem 1.2 follows Federer's strategy in [19] that a basic Poincaré inequality upgrades to the stronger Boxing inequality via a covering argument. As it is straightforward to get some form of fractional Poincaré inequality from the definition of the Gagliardo seminorm, this is sufficient to obtain an elementary proof of a weaker variant of our result, one without the proper asymptotics (see also Theorem 2.1 in [35], with $\theta = 1$ and $l = \alpha$). However, a second aspect of our argument is to connect the asymptotics as $\alpha$ tends to $0$ and $1$. In the former regime we rely on the behavior at infinity of the function $1/|x|^{d+\alpha}$, while in the latter we utilize a Poincaré inequality for the fractional perimeter with the right behavior of the constant as $\alpha$ tends to $1$. This exploits the fact that after one has established an inequality for each $\alpha$ fixed, the asymptotics can be handled by a $\Gamma$-convergence counterpart of the second limit of (1.3). This approach bears some analogy with the $\Gamma$-convergence of non-local functionals on BV-functions from [38, Corollary 8]; see also [7].

Theorems 1.1 and 1.2 have powerful implications in the study of Sobolev functions which can be understood through a functional formulation of the Boxing inequality. To this end we recall the definition of the Hausdorff content of dimension $d - \alpha$, which for any set $A \subset \mathbb{R}^d$ is given by

$$(1.4) \qquad \mathcal{H}^{d-\alpha}_\infty(A) := \inf\left\{\sum_{i=0}^\infty \omega_{d-\alpha} r_i^{d-\alpha} : A \subset \bigcup_{i=0}^\infty B_{r_i}(x_i)\right\},$$

where $\omega_{d-\alpha} := \pi^{(d-\alpha)/2}/\Gamma\left(\frac{d-\alpha}{2} + 1\right)$. Then integration of a function $u : \mathbb{R}^d \to \mathbb{R}$ with respect to the Hausdorff content defines the Choquet integral as

$$\int_{\mathbb{R}^d} |u| \, \mathrm{d}\mathcal{H}^{d-\alpha}_\infty := \int_0^\infty \mathcal{H}^{d-\alpha}_\infty(\{|u| > t\}) \, \mathrm{d}t.$$

With these tools we are ready to state our

**Theorem 1.3.** *Let $\alpha \in (0, 1)$. There exists a constant $C = C(d) > 0$ such that*

$$\int_{\mathbb{R}^d} |\varphi| \, \mathrm{d}\mathcal{H}^{d-\alpha}_\infty \leq C\alpha(1-\alpha)[\varphi]_{W^{\alpha,1}(\mathbb{R}^d)},$$

*for every $\varphi \in C_c^\infty(\mathbb{R}^d)$.*

**Remark 1.4.** The estimate in Theorem 1.3 extends to all $u \in W^{\alpha,1}(\mathbb{R}^d)$, though one should rely on the precise representative $u^*$ in the Choquet integral; see Section 3 below.

The proof of Theorem 1.3 follows the simple yet beautiful idea of Federer and Fleming [21] and Maz'ya [29] that the combination of the coarea formula and the isoperimetric inequality yield Sobolev's inequality in the $L^1$ regime. The coarea formula for $W^{\alpha,1}(\mathbb{R}^d)$ is given in (1.2), while the Boxing inequality provides a replacement of the isoperimetric component via

the estimate

$$\mathcal{H}^{d-\alpha}_{\infty}(U) \le \omega_{d-\alpha} \, C\alpha(1-\alpha) P_\alpha(U), \tag{1.5}$$

which is an easy consequence of Theorem 1.2.

Theorem 1.3 is a strong form of differential inequality for functions $u \in W^{\alpha,1}(\mathbb{R}^d)$ that captures in a precise way the fine properties of $u$. From it one can deduce a variety of integral estimates, as well as quantify the size of the Lebesgue set of $u$. For example, it firstly implies the fractional Sobolev embedding on the scale of Lebesgue $L^p$ spaces, which is

**Corollary 1.5.** *Let $\alpha \in (0,1)$. There exists a constant $C = C(d) > 0$ such that*

$$\|u\|_{L^{\frac{d}{d-\alpha}}(\mathbb{R}^d)} \le C\alpha(1-\alpha)[u]_{W^{\alpha,1}(\mathbb{R}^d)}$$

*for every $u \in W^{\alpha,1}(\mathbb{R}^d)$.*

That is, in place of interpolation [28], rearrangements to obtain an isoperimetric inequality [5, 25], or various other methods [12, 32], one obtains directly a Sobolev inequality that is stable in the limit as $\alpha$ tends to $0$ and $1$. More than this, Theorem 1.3 encodes a trace inequality that enables one to control the integral of functions $u \in W^{\alpha,1}(\mathbb{R}^d)$ along lower dimension objects. In practice this takes the form of the a priori inequality

$$\int_{\mathbb{R}^d} |\varphi| \, \mathrm{d}\mu \le C\alpha(1-\alpha)[\varphi]_{W^{\alpha,1}(\mathbb{R}^d)} \tag{1.6}$$

for every $\varphi \in C_c^\infty(\mathbb{R}^d)$ and every nonnegative Borel measure $\mu$ in $\mathbb{R}^d$ that satisfies

$$\mu(B_r(x)) \le \omega_{d-\alpha} r^{d-\alpha} \tag{1.7}$$

for all balls $B_r(x) \subset \mathbb{R}^d$. This estimate is a direct consequence of (1.5), with the same constant, since (1.7) yields the comparison $\mu \le \mathcal{H}^{d-\alpha}_\infty$, while again inequality (1.6) extends to functions in $W^{\alpha,1}(\mathbb{R}^d)$ when one utilizes precise representatives. In this form, one can readily prove a number of other inequalities for Sobolev functions. For example, one observes that $\mu = g \, \mathrm{d}x$, with $g$ in the Lorentz space $L^{\frac{d}{\alpha},\infty}(\mathbb{R}^d)$ of weak $L^{\frac{d}{\alpha}}$ functions, or $\mu = \mathrm{d}x/|x|^\alpha$ both satisfy such a ball-growth condition. For the former, one obtains by duality a result stronger than Corollary 1.5, namely the embedding of $W^{\alpha,1}(\mathbb{R}^d)$ into the Lorentz space $L^{\frac{d}{d-\alpha},1}(\mathbb{R}^d)$ (see [6]), while the latter yields Hardy's inequality. We refer the reader to [33] for further discussions about trace inequalities.

The analogue to Theorem 1.3 and (1.6) on bounded domains for $\alpha = 1$ is due to Meyers and Ziemer in their paper [34] on Poincaré-Wirtinger inequalities. We can pursue such inequalities in this framework, which require a modification of our argument to obtain an estimate similar to (1.5) involving the relative perimeter $P_\alpha(U, \Omega)$. The result of this analysis is the fractional Sobolev-Poincaré inequality given by our

**Theorem 1.6.** *Let $\alpha \in (0, 1)$ and $\Omega \subset \mathbb{R}^d$ be a smooth connected bounded open set. There exists a constant $C = C(d, \Omega) > 0$ such that*

$$\left\|\varphi - \fint_\Omega \varphi\right\|_{L^1(\Omega, \mathrm{d}\mu)} \leq C(1-\alpha) \int_\Omega \int_\Omega \frac{|\varphi(y) - \varphi(z)|}{|y-z|^{\alpha+d}} \, \mathrm{d}y \, \mathrm{d}z,$$

*for every $\varphi \in C^\infty(\overline{\Omega})$ and every nonnegative Borel measure $\mu \leq \mathcal{H}_\infty^{d-\alpha}$ in $\Omega$.*

The absence of the factor of $\alpha$ in the preceding estimate is expected since the double integral in the right-hand side converges to a finite limit as $\alpha$ tends to 0. In contrast, its appearance in Theorem 1.3 comes from the behavior of the function $1/|x|^{\alpha+d}$ for $|x|$ large, while now $\Omega$ is bounded.

Further applications of the functional formulation of the Boxing inequality given in Theorem 1.3 can be found when studying the fine properties of functions in the Sobolev space $W^{\alpha,1}(\mathbb{R}^d)$. We explore this in Sections 2 and 3, where we are additionally interested in the larger range $\alpha \in [0, d]$. In this interval we have the following extension of Theorem 1.3, which is our

**Theorem 1.7.** *Let $\alpha \in [0, d]$. There exists a constant $C = C(d, \alpha) > 0$ such that*

$$\int_{\mathbb{R}^d} |\varphi| \, \mathrm{d}\mathcal{H}_\infty^{d-\alpha} \leq C \, [\nabla^k \varphi]_{W^{\alpha-k,1}(\mathbb{R}^d)},$$

*for every $\varphi \in C_c^\infty(\mathbb{R}^d)$. Here $k = \lfloor \alpha \rfloor$ is the integer part of $\alpha$, with the convention that $[\nabla^k \varphi]_{W^{0,1}(\mathbb{R}^d)} := \|\nabla^k \varphi\|_{L^1(\mathbb{R}^d)}$. Moreover, for $\alpha \in (0, d) \setminus \mathbb{N}$, the constant above satisfies*

(1.8) $$C \leq C'(d) \operatorname{dist}(\alpha, \mathbb{N}).$$

A few comments regarding Theorem 1.7. First, as in the proof of Theorem 1.2, for $\alpha \in (k, k+1)$ we handle the cases $\alpha \to k^+$ and $\alpha \to (k+1)^-$ by different arguments. In the case of the former, we rely on the existence of bounded solutions of the divergence equation

$$-\operatorname{div} Y = \nu$$

for nonnegative measures $\nu$ in $\mathbb{R}^d$ which satisfy the ball-growth condition $\nu(B_r(x)) \leq r^{d-1}$. This divergence equation is dual to inequality (1.6) for $\alpha = 1$, and so we do not need to invoke any further results on special solutions to such an equation (see for example [10, 13]). The latter case $\alpha \to (k+1)^-$ is a direct consequence of lifting Theorem 1.3 via the mapping properties of Riesz potentials (see Lemma 4.6 below) and inherits the correct scaling from our Theorem 1.3. Therefore the ability to handle decreasing limits for $\mathcal{H}_\infty^{d-\alpha}$, our estimate, and the asymptotics of the Gagliardo semi-norms (see [32]) imply by continuity D. Adams' estimate

(1.9) $$\int_{\mathbb{R}^d} |\varphi| \, \mathrm{d}\mathcal{H}_\infty^{d-k} \leq C \|\nabla^k \varphi\|_{L^1(\mathbb{R}^d)}$$

for all integer orders $k \in \{1, \ldots, d\}$; see [2].

Second, the key point of our estimate is found in (1.8). If one is willing to accept $H^1$-BMO duality and dispense with this bound, then one can deduce such a result from Adams' capacitary inequalities in [2]; we explore this connection in Section 5. Finally, the Sobolev inequality in Theorem 1.7 extends the range of Corollary 1.5, that one has

$$(1.10) \qquad \|u\|_{L^{\frac{d}{d-\alpha}}(\mathbb{R}^d)} \leq C'(d) \operatorname{dist}(\alpha, \mathbb{N}) \, [\nabla^k u]_{W^{\alpha-k,1}(\mathbb{R}^d)},$$

for every $\alpha \in (0, d) \setminus \mathbb{N}$ and $u \in W^{\alpha,1}(\mathbb{R}^d)$, and where again $k = \lfloor \alpha \rfloor$.

The plan of the paper is as follows. In Section 2 we explore the relationship of Hausdorff content and a capacity that is intrinsically associated to $W^{\alpha,1}(\mathbb{R}^d)$. In Section 3 we show the implications of our results on the fine properties of functions $W^{\alpha,1}(\mathbb{R}^d)$. In Section 4 we prove Theorems 1.2, 1.3, 1.6 and 1.7 and Corollary 1.5 above. Finally, in Section 5 we relate our results with those obtained by Adams in [2]. In particular we show how $H^1$-BMO duality can be used to deduce Theorem 1.7, however missing the asymptotics we obtain in our results.

## 2. The Equivalence of Capacity and Content

For every $\alpha \geq 0$ we define the $(\alpha, 1)$-capacity of a compact set $K \subset \mathbb{R}^d$ by

$$\operatorname{Cap}_{\alpha,1}(K) := \inf \Big\{ [\nabla^k \varphi]_{W^{\theta,1}(\mathbb{R}^d)} : \varphi \in C_c^\infty(\mathbb{R}^d), \ \varphi \geq 0 \text{ in } \mathbb{R}^d, \ \varphi > 1 \text{ on } K \Big\},$$

where

$$k = \lfloor \alpha \rfloor \in \mathbb{N} \qquad \text{and} \qquad \theta = \alpha - k \in [0, 1).$$

We take the convention that $\nabla^0 \varphi = \varphi$ and $[\nabla^k \varphi]_{W^{0,1}(\mathbb{R}^d)} = \|\nabla^k \varphi\|_{L^1(\mathbb{R}^d)}$.

One then extends the capacity by regularity as follows: For an open set $U \subset \mathbb{R}^d$, the capacity of $U$ is the supremum of capacities of compact subsets $K \subset U$:

$$\operatorname{Cap}_{\alpha,1}(U) := \sup_{K \subset U} \operatorname{Cap}_{\alpha,1}(K).$$

Then for a general set $A \subset \mathbb{R}^d$ the capacity is the infimum of the capacities of open supersets $U \supset A$,

$$\operatorname{Cap}_{\alpha,1}(A) := \inf_{U \supset A} \operatorname{Cap}_{\alpha,1}(U).$$

The only point is to check that the definition on compact sets is consistent, but this is indeed the case (see for example Proposition A.6 in [40]).

Capacities are outer measures that arise commonly in the study of partial differential equations when one seeks to quantify the size of the set for which certain fine properties hold. For instance, for a function $u$ in the Sobolev space $W^{\alpha,1}(\mathbb{R}^d)$, the set of points $x \in \mathbb{R}^d$ for which the limit

$$\lim_{r \to 0} \fint_{B_r(x)} u \tag{2.1}$$

fails to exist has $(\alpha, 1)$-capacity zero; see Section 3 below. When $\alpha = 0$ this statement is none other than the classical Lebesgue Differentiation Theorem, which states that for every function $u \in L^1(\mathbb{R}^d)$ the limit (2.1) exists except on a subset with Lebesgue measure zero. Indeed, the $(0, 1)$-capacity yields an alternative construction of the Lebesgue outer measure.

For $\alpha > 0$ one expects that the better regularity of a $W^{\alpha,1}$-function yield a smaller exceptional set with regard to the existence of the limit (2.1). For example, functions in the Sobolev space $W^{\alpha,1}(\mathbb{R}^d)$ for $\alpha \geq d$ are continuous, so that the limit (2.1) holds everywhere. This is confirmed by the quantification by capacity, since every nonempty set has positive capacity in this range of $\alpha$. However, when $\alpha \in (0, d)$ it is not obvious from the definition of the capacity how small this set is, as it is intrinsically given by a Sobolev semi-norm. It is therefore desirable to connect these capacities for $\alpha \in (0, d)$ with a geometric description of size of sets in $\mathbb{R}^d$.

Some natural geometric objects for capturing the size of sets in this respect are the Hausdorff measures, which give a way of assigning a notion of lower-order dimension to sets of Lebesgue measure zero. The first step in their construction is the definition of certain outer measures defined for example using coverings with balls, which are measured as if they were $s$-dimensional objects: For any set $A \subset \mathbb{R}^d$, define

$$\mathcal{H}^s_\delta(A) := \inf\left\{\sum_{n=0}^\infty \omega_s r_n^s : A \subset \bigcup_{n=0}^\infty B_{r_n}(x_n), r_n \leq \delta\right\}. \tag{2.2}$$

From here one can define the $s$-dimensional Hausdorff measure of $A$ as the non-decreasing limit

$$\mathcal{H}^s(A) = \lim_{\delta \to 0} \mathcal{H}^s_\delta(A).$$

The Hausdorff measure $\mathcal{H}^s$ provides a geometric object for performing integration on an embedded surface of dimension $s \in \mathbb{N}$. For example, an application of this construction to the dimension $s = d - 1$ and restricting to the surface of the unit sphere $\partial B_1$, one obtains the standard surface measure on the sphere. More generally, the Hausdorff measures provide a way of assigning a measure to sets between the integer dimensions, like the Cantor middle-third set, which has finite but non-zero $\log(2)/\log(3)$-dimensional Hausdorff measure.

While the Hausdorff measures describe the size of and assign measure to sets in $\mathbb{R}^d$ aptly, the $\mathcal{H}^s$-measure of a non-empty ball in $\mathbb{R}^d$ is infinite for every $s \in [0, d)$. Thus they are too large to make setwise comparison with capacities. A smaller object than the Hausdorff measures, and also more suitable to comparison with the homogeneous $(\alpha, 1)$-capacity, is the Hausdorff content $\mathcal{H}_\infty^{d-\alpha}$ as defined by (1.4) (which is (2.2) with $\delta = \infty$ and $s = d - \alpha$). Indeed,

**Theorem 2.1.** *For every $\alpha \in (0, d)$ and every $A \subset \mathbb{R}^N$,*
$$\operatorname{Cap}_{\alpha,1}(A) \sim \mathcal{H}_\infty^{d-\alpha}(A).$$

The case $\alpha = 1$ was established by Meyers and Ziemer [34] and the cases $\alpha \in \{2, \ldots, d-1\}$ by Adams [2], based on an induction argument by Maz'ya [30]. We prove Theorem 2.1 in Section 4 below.

One can define $(\alpha, p)$-capacities for any $p > 1$ in analogy with the definition of $\operatorname{Cap}_{\alpha,1}$ by minimization of the Gagliardo energy
$$[\nabla^k \varphi]_{W^{\alpha-k,p}(\mathbb{R}^d)}^p := \int_{\mathbb{R}^d} \int_{\mathbb{R}^d} \frac{|\nabla^k \varphi(y) - \nabla^k \varphi(z)|^p}{|y-z|^{(\alpha-k)p+d}} \, dy \, dz$$

for $\alpha - k \in (0, 1)$. Nevertheless, an equivalence between capacity and content as in Theorem 2.1 holds exclusively for $p = 1$; see Chapter 5 in [4]. One direction of the comparison, which has a true counterpart for $p > 1$, is an easy consequence of having the associated scaling. Indeed, one has
$$\operatorname{Cap}_{\alpha,1}(B_r(x)) = P_\alpha(B_1) \, r^{d-\alpha}$$

when $\alpha \in (0, 1)$, and similarly $\operatorname{Cap}_{1,1}(B_r(x)) = \operatorname{Per}(B_1) \, r^{d-1}$; see e.g. [44, 45]. An analogous behavior is also valid for any $\alpha \in (0, d)$ and implies

(2.3) $$\operatorname{Cap}_{\alpha,1}(A) \leq C \mathcal{H}_\infty^{d-\alpha}(A),$$

for every $A \subset \mathbb{R}^d$, with constant $C = \operatorname{Cap}_{\alpha,1}(B_1)$. The reverse comparison is more subtle, and relies strongly on the $L^1$-nature of the Boxing inequality; see the end of Section 4.

With these ingredients, we can return to the question of the fine properties of a function in the Sobolev space $W^{\alpha,1}(\mathbb{R}^d)$. Theorem 2.1 allows one to quantify the size of the exceptional set of a Sobolev function in terms of the Hausdorff measure (and not the content). This follows from the remarkable but easy fact concerning the Hausdorff measure $\mathcal{H}^{d-\alpha}$ and the Hausdorff content $\mathcal{H}_\infty^{d-\alpha}$ that they have the same negligible sets. Therefore a result concerning the equivalence of the $(\alpha, 1)$-capacity and the $\mathcal{H}_\infty^{d-\alpha}$-content is sufficient to guarantee that sets of $(\alpha, 1)$-capacity zero are $\mathcal{H}^{d-\alpha}$-negligible. In particular, we have that

(2.4) $$\mathcal{H}^{d-\alpha}(A) = 0 \iff \operatorname{Cap}_{\alpha,1}(A) = 0,$$

which implies that for any $u \in W^{\alpha,1}(\mathbb{R}^d)$ the set of points where the limit (2.1) does not exist is negligible with respect to $\mathcal{H}^{d-\alpha}$. The equivalence (2.4) was first observed for the $(1, 1)$-capacity by Fleming in his pioneer paper [24].

## 3. Fine Properties of Sobolev Functions

We explain in this section why the exceptional set of a function $u$ in the Sobolev space $W^{\alpha,1}(\mathbb{R}^d)$ has $(\alpha, 1)$-capacity zero. This property is proved in a standard way once a suitable maximal-function inequality is established; see (3.4). We use such an information on $u$ to extend Theorem 1.7, which includes Corollary 1.5 and estimate (1.6), to any function in $W^{\alpha,1}(\mathbb{R}^d)$.

We first recall that $x \in \mathbb{R}^d$ is a Lebesgue point of a function $u \in L^1(\mathbb{R}^d)$ if there exists $a \in \mathbb{R}$ such that

$$\text{(3.1)} \qquad \lim_{r \to 0} \fint_{B_r(x)} |u - a| = 0.$$

We denote by $\mathcal{L}_u$ the set of Lebesgue points of $u$. Assigning the value of $a$ to $x \in \mathcal{L}_u$ yields a function $u^* : \mathcal{L}_u \to \mathbb{R}$, $u^*(x) := a$, which is called the precise representative of $u$. We thus have

$$\text{(3.2)} \qquad \lim_{r \to 0} \fint_{B_r(x)} |u - u^*(x)| = 0$$

and then

$$\text{(3.3)} \qquad \lim_{r \to 0} \fint_{B_r(x)} u = u^*(x).$$

An advantage of having (3.2) over (3.3) is that other averaging processes yield the same value $u^*(x)$. For example, given any $\rho \in C_c^\infty(\mathbb{R}^d)$ such that $\int_{\mathbb{R}^d} \rho = 1$, it follows from (3.2) that

$$\lim_{r \to 0} \frac{1}{r^d} \int_{\mathbb{R}^d} \rho\Big(\frac{x-y}{r}\Big) u(y) \, \mathrm{d}y = u^*(x).$$

By the Lebesgue Differentiation Theorem,

$$\mathrm{Cap}_{0,1}\,(\mathbb{R}^d \setminus \mathcal{L}_u) = |\mathbb{R}^d \setminus \mathcal{L}_u| = 0$$

and $u^* = u$ almost everywhere in $\mathbb{R}^d$. We are interested in estimating the size of the exceptional set $\mathbb{R}^d \setminus \mathcal{L}_u$ in terms of the capacity when $u$ is a Sobolev function:

**Proposition 3.1.** *For every $\alpha \in (0, d)$ and every $u \in W^{\alpha,1}(\mathbb{R}^d)$,*

$$\mathrm{Cap}_{\alpha,1}\,(\mathbb{R}^d \setminus \mathcal{L}_u) = \mathcal{H}_\infty^{d-\alpha}(\mathbb{R}^d \setminus \mathcal{L}_u) = 0.$$

The main ingredient follows from the strong-type estimate for the maximal function:

**Lemma 3.2.** *For every $\alpha \in (0, d)$ and every $u \in W^{\alpha,1}(\mathbb{R}^d)$, one has*

$$\int_{\mathbb{R}^d} \mathcal{M}u \, \mathrm{d}\mathcal{H}_\infty^{d-\alpha} \le C[u]_{W^{\alpha,1}(\mathbb{R}^d)},$$

*where $\mathcal{M}u : \mathbb{R}^d \to [0, \infty]$ is the Hardy-Littlewood maximal function defined by*

$$\mathcal{M}u(x) := \sup_{r > 0} \fint_{B_r(x)} |u|.$$

Being a supremum of a family of continuous functions, the maximal function is lower semicontinuous, hence all the sets $\{\mathcal{M}u > t\}$ are open. A consequence of this lemma is the weak-type estimate

$$(3.4) \qquad \mathcal{H}_\infty^{d-\alpha}(\{\mathcal{M}u > t\}) \leq \frac{C}{t}[u]_{W^{\alpha,1}(\mathbb{R}^d)},$$

for every $t > 0$. Notice that for all $x \in \mathcal{L}_u$,

$$|u^*(x)| \leq \liminf_{r \to 0} \fint_{B_r(x)} |u| \leq \mathcal{M}u(x).$$

Since $\mathcal{H}_\infty^{d-\alpha}(\mathbb{R}^d \setminus \mathcal{L}_u) = 0$, applying Lemma 3.2 we can express the inequality in Theorem 1.7 for every $u \in W^{\alpha,1}(\mathbb{R}^d)$ as

$$(3.5) \qquad \int_{\mathbb{R}^d} |u^*| \, d\mathcal{H}_\infty^{d-\alpha} \leq \int_{\mathbb{R}^d} \mathcal{M}u \, d\mathcal{H}_\infty^{d-\alpha} \leq C \, [u]_{W^{\alpha,1}(\mathbb{R}^d)},$$

where $C > 0$ also satisfies (1.8) for $\alpha \in (0, d) \setminus \mathbb{N}$. Further applications of these estimates and the equivalence of $\text{Cap}_{\alpha,1}$ and $\mathcal{H}_\infty^{d-\alpha}$ in the study of properties of $W^{\alpha,1}$-functions will be investigated in a forthcoming work.

*Proof of Lemma 3.2.* By the functional Boxing inequality (Theorem 1.7), we can write

$$(3.6) \qquad \int_{\mathbb{R}^d} |\varphi| \, d\mathcal{H}_\infty^{d-\alpha} \leq C_1 [\varphi]_{W^{\alpha,1}(\mathbb{R}^d)},$$

for every $\varphi \in C_c^\infty(\mathbb{R}^d)$. We also rely on Adams' maximal-function estimate involving the Choquet integral with respect to $\mathcal{H}_\infty^{d-\alpha}$:

$$(3.7) \qquad \int_{\mathbb{R}^d} \mathcal{M}\varphi \, d\mathcal{H}_\infty^{d-\alpha} \leq C_2 \int_{\mathbb{R}^d} |\varphi| \, d\mathcal{H}_\infty^{d-\alpha};$$

an elegant proof of this estimate is due to Orobitg and Verdera [36]. Combining (3.6) and (3.7), we have

$$(3.8) \qquad \int_{\mathbb{R}^d} \mathcal{M}\varphi \, d\mathcal{H}_\infty^{d-\alpha} \leq C_3 [\varphi]_{W^{\alpha,1}(\mathbb{R}^d)}.$$

Given $u \in W^{\alpha,1}(\mathbb{R}^d)$, we apply this estimate to a sequence $(\varphi_j)_{j \in \mathbb{N}}$ in $C_c^\infty(\mathbb{R}^d)$ which converges to $u$ in $W^{\alpha,1}(\mathbb{R}^d)$. Some care is needed to justify the limit in the left-hand side: One can proceed along the lines of the proof of Lemma 9.8 in [40], but the argument there relies on the strong subadditivity of Choquet's capacity. It seems unlikely that the spherical Hausdorff content $\mathcal{H}_\infty^{d-\alpha}$ is strongly subadditive, so instead one uses the dyadic Hausdorff content $\widehat{\mathcal{H}}_\infty^{d-\alpha}$ defined for every $A \subset \mathbb{R}^d$ by

$$(3.9) \qquad \widehat{\mathcal{H}}_\infty^{d-\alpha}(A) = \inf \left\{ \sum_{i=0}^\infty \ell_i^{d-\alpha} : A \subset \text{int} \bigcup_{i=0}^\infty Q_i \right\},$$

where the infimum is computed over all sequences of closed dyadic cubes $Q_i$, and $\ell_i$ denotes the side length of the cube $Q_i$. One observes that $\mathcal{H}_\infty^{d-\alpha} \sim$ $\widehat{\mathcal{H}}_\infty^{d-\alpha}$ and $\widehat{\mathcal{H}}_\infty^{d-\alpha}$ is strongly subadditive in the sense that

$$\widehat{\mathcal{H}}_\infty^{d-\alpha}(A \cup B) + \widehat{\mathcal{H}}_\infty^{d-\alpha}(A \cap B) \le \widehat{\mathcal{H}}_\infty^{d-\alpha}(A) + \widehat{\mathcal{H}}_\infty^{d-\alpha}(B),$$

for every $A, B \subset \mathbb{R}^d$; see [46]. Thus, $\widehat{\mathcal{H}}_\infty^{d-\alpha}$ satisfies the conclusion of the Increasing Set Lemma and we have

$$(3.10) \qquad \lim_{k \to \infty} \widehat{\mathcal{H}}_\infty^{d-\alpha}(A_k) = \widehat{\mathcal{H}}_\infty^{d-\alpha}\Big(\bigcup_{k=0}^\infty A_k\Big),$$

for every non-decreasing sequence $(A_k)_{k \in \mathbb{N}}$ of subsets of $\mathbb{R}^d$. By the argument in the proof of Lemma 9.8 in [40] we then deduce that

$$\widehat{\mathcal{H}}_\infty^{d-\alpha}(\{\mathcal{M}u > t\}) \le \liminf_{j \to \infty} \widehat{\mathcal{H}}_\infty^{d-\alpha}(\{\mathcal{M}\varphi_j > t\}),$$

for every $t > 0$. By Fatou's lemma we get

$$\int_0^\infty \widehat{\mathcal{H}}_\infty^{d-\alpha}(\{\mathcal{M}u > t\})\,\mathrm{d}t \le \liminf_{j \to \infty} \int_0^\infty \widehat{\mathcal{H}}_\infty^{d-\alpha}(\{\mathcal{M}\varphi_j > t\})\,\mathrm{d}t.$$

In view of (3.8) and the equivalence between the Hausdorff contents $\mathcal{H}_\infty^{d-\alpha}$ and $\widehat{\mathcal{H}}_\infty^{d-\alpha}$, the conclusion follows. $\square$

*Proof of Proposition 3.1.* The existence of the limit in (3.1) at a point $x \in \mathbb{R}^d$ is equivalent to the Cauchy condition:

$$\lim_{(r,s) \to (0,0)} \fint_{B_r(x)} \fint_{B_s(x)} |u(y) - u(z)|\,\mathrm{d}y\,\mathrm{d}z = 0;$$

see [40, Lemma 8.8]. We may then assert that $\mathbb{R}^d \setminus \mathcal{L}_u = \bigcup_{\lambda > 0} A_\lambda$, where

$$A_\lambda = \Big\{ x \in \mathbb{R}^d : \limsup_{(r,s) \to (0,0)} \fint_{B_r(x)} \fint_{B_s(x)} |u(y) - u(z)|\,\mathrm{d}y\,\mathrm{d}z > \lambda \Big\}.$$

One can now proceed along the lines of the proof of Proposition 8.6 in [40] using the weak-type estimate (3.4) for the maximal function and the density of $C_c^\infty(\mathbb{R}^d)$ in $W^{\alpha,1}(\mathbb{R}^d)$ to deduce that $\mathcal{H}_\infty^{d-\alpha}(A_\lambda) = 0$, for every $\lambda > 0$. The monotonicity of the family $(A_\lambda)_{\lambda > 0}$ allows one to write

$$\bigcup_{\lambda > 0} A_\lambda = \bigcup_{n=0}^\infty A_{\lambda_n},$$

where $(\lambda_n)_{n \in \mathbb{N}}$ is any sequence of positive numbers that converges to $0$. The countable subadditivity of $\mathcal{H}_\infty^{d-\alpha}$ and estimate (2.3) then imply that

$$\mathrm{Cap}_{\alpha,1}(\mathbb{R}^d \setminus \mathcal{L}_u) \le C \mathcal{H}_\infty^{d-\alpha}(\mathbb{R}^d \setminus \mathcal{L}_u) = 0. \qquad \square$$

4. Proofs of the Main Results

Before proving Theorem 1.2, we first show that a relative isoperimetric inequality holds uniformly with respect to $\alpha \in (0,1)$. We begin with the behavior near $\alpha = 1$:

**Lemma 4.1.** *Given $\gamma \in (0,1)$, there exists a constant $C' = C'(d, \gamma) > 0$ such that*
$$r^{d-\alpha} \leq C'(1-\alpha) \int_A \int_{B_r \setminus A} \frac{\mathrm{d}y \, \mathrm{d}z}{|y-z|^{\alpha+d}}$$
*for every $r > 0$, every $\alpha \in (0,1)$, and every Borel set $A \subset B_r$ such that $|A|/|B_r| = \gamma$.*

*Proof.* By a scaling argument, it suffices to consider the case where $r = 1$. The estimate follows from the fractional Poincaré inequality

(4.1) $$\int_{B_1} \int_{B_1} |u(y) - u(z)| \, \mathrm{d}y \, \mathrm{d}z \leq C_1 (1-\alpha) [u]_{W^{\alpha,1}(B_1)}$$

for every $\alpha \in (0,1)$ and $u \in L^1(B_1)$. Indeed, it suffices to apply this inequality with $u = \chi_A$ to get
$$2|A||B_1 \setminus A| \leq C_1(1-\alpha)[\chi_A]_{W^{\alpha,1}(B_1)},$$
and then we get the estimate we want since $|A|/|B_1| = \gamma$ and $r = 1$:
$$|B_1|^2 \gamma(1-\gamma) \leq C_1(1-\alpha) \int_A \int_{B_1 \setminus A} \frac{\mathrm{d}y \, \mathrm{d}z}{|y-z|^{\alpha+d}}.$$

I Inequality (4.1), without the coefficient $(1-\alpha)$, is a straightforward consequence of the definition of the Gagliardo seminorm. We can thus focus in the range $\alpha \in [1/2, 1)$, and then (4.1) follows from the inequality

(4.2) $$\int_{B_1} \left| u(x) - \fint_{B_1} u \right| \mathrm{d}x \leq C_2(1-\alpha)[u]_{W^{\alpha,1}(B_1)},$$

due to Bourgain, Brezis and Mironescu [12, p. 80]. The proof in [12] is based on a Paley-Littlewood decomposition of $u$. A more elementary approach is presented in [39, Section 6] and relies on a standard contradiction argument based on the relative compactness in $L^1(B_1)$ of a sequence of functions $(u_n)_{n \in \mathbb{N}}$ that satisfies
$$(1-\alpha_n)[u_n]_{W^{\alpha_n,1}(B_1)} \leq C_2,$$
where $(\alpha_n)_{n \in \mathbb{N}}$ is a sequence in $(0,1)$ that converges to 1; see [11, Theorem 4] and also [39, Theorems 1.2 and 1.3] for a generalization. Observe that by the triangle inequality one has
$$\int_{B_1} \int_{B_1} |u(y) - u(z)| \, \mathrm{d}x \, \mathrm{d}y \leq 2|B_1| \int_{B_1} \left| u(x) - \fint_{B_1} u \right| \mathrm{d}x.$$

Combining this estimate with (4.2), one gets (4.1) for $\alpha \in [1/2, 1)$ and this completes the proof of the lemma. □

We next focus on the relative isoperimetric inequality near $\alpha = 0$. The analysis in this case relies on the behavior of the potential $1/|x|^{\alpha+d}$ for $|x|$ large.

**Lemma 4.2.** *Given $\gamma \in (0,1)$ and a bounded Borel set $E \subset \mathbb{R}^d$, there exists a constant $C'' = C''(d, \gamma) > 0$ such that*
$$r^{d-\alpha} \le C'' \alpha \int_{B_r \cap E} \int_{\mathbb{R}^d \setminus E} \frac{\mathrm{d}y\,\mathrm{d}z}{|y-z|^{\alpha+d}},$$
*for every $\alpha \in (0,1)$, where $r > 0$ is such that*
$$|B_r \cap E|/|B_r| = \gamma \quad \text{and} \quad |B_s \cap E|/|B_s| \le \gamma \quad \text{for every } s \ge r.$$

*Proof.* Denoting by $I$ the double integral in the right-hand side, we have
$$I \ge \int_{B_r \cap E} \left( \int_{(\mathbb{R}^d \setminus E) \setminus B_{2r}} \frac{\mathrm{d}z}{|y-z|^{\alpha+d}} \right) \mathrm{d}y.$$

For every $y \in B_r$ and $z \in \mathbb{R}^d \setminus B_{2r}$, by the triangle inequality we have
$$|y - z| \le |y| + |z| \le r + |z| \le \frac{3}{2}|z|.$$

Thus, for every $\alpha \in (0,1)$ we have
$$I \ge \left(\frac{2}{3}\right)^{1+d} \int_{B_r \cap E} \left( \int_{(\mathbb{R}^d \setminus E) \setminus B_{2r}} \frac{\mathrm{d}z}{|z|^{\alpha+d}} \right) \mathrm{d}y$$
$$= \left(\frac{2}{3}\right)^{1+d} |B_r \cap E| \int_{(\mathbb{R}^d \setminus E) \setminus B_{2r}} \frac{\mathrm{d}z}{|z|^{\alpha+d}} = C_1 \gamma\, r^d \int_{(\mathbb{R}^d \setminus E) \setminus B_{2r}} \frac{\mathrm{d}z}{|z|^{\alpha+d}}.$$

To conclude the proof, it suffices to show that
$$\int_{(\mathbb{R}^d \setminus E) \setminus B_{2r}} \frac{\mathrm{d}z}{|z|^{\alpha+d}} \ge \frac{C_2}{\alpha\, r^\alpha}.$$

For this purpose, using Cavalieri's principle we first rewrite
$$\int_{(\mathbb{R}^d \setminus E) \setminus B_{2r}} \frac{\mathrm{d}z}{|z|^{\alpha+d}} = (\alpha + d) \int_{2r}^\infty \frac{|(B_s \setminus E) \setminus B_{2r}|}{s^{\alpha+d}} \frac{\mathrm{d}s}{s}.$$

Using the assumption on $r$, we can then estimate
$$|(B_s \setminus E) \setminus B_{2r}| \ge |B_s \setminus E| - |B_{2r}| \ge (1-\gamma)|B_s| - |B_{2r}| = \omega_d \left[ (1-\gamma) - \left(\frac{2r}{s}\right)^d \right] s^d.$$

Take a fixed number $\lambda > 2$ such that $(1-\gamma) > (2/\lambda)^d$. It thus follows that for every $s \ge \lambda r$ we have
$$|(B_s \setminus E) \setminus B_{2r}| \ge C_3 s^d.$$

We finally get
$$\int_{(\mathbb{R}^d \setminus E) \setminus B_{2r}} \frac{\mathrm{d}z}{|z|^{\alpha+d}} \ge d \int_{\lambda r}^\infty \frac{|(B_s \setminus E) \setminus B_{2r}|}{s^{\alpha+d}} \frac{\mathrm{d}s}{s} \ge d \int_{\lambda r}^\infty \frac{C_3 s^d}{s^{\alpha+d}} \frac{\mathrm{d}s}{s} = \frac{C_4}{\alpha\, r^\alpha}.$$

This concludes the proof of the lemma. $\square$

*Proof of Theorem 1.2.* Given $x \in U$, let $r = r(x) > 0$ be the largest number such that
$$\frac{|B_r(x) \cap U|}{|B_r(x)|} = \frac{1}{2}.$$
The existence of an $r$ for which the equality holds follows from the Intermediate Value Theorem: The quantity in the left-hand side is continuous with respect to $r > 0$, equals $1$ for $r$ small since $x \in U$ and $U$ is open, and converges to $0$ as $r \to \infty$ since $U$ is bounded. A continuity argument shows that a largest solution $r$ indeed exists. Observe that such a choice of $r$ does not depend on $\alpha$.

Applying Lemma 4.1 on the ball $B_r(x)$ to the open set $A = B_r(x) \cap U$, we have
$$r^{d-\alpha} \leq C_1(1-\alpha) \int_{B_r(x) \cap U} \int_{\mathbb{R}^d \setminus U} \frac{\mathrm{d}y\,\mathrm{d}z}{|y-z|^{\alpha+d}}.$$

By Lemma 4.2 with $E = U$, we also have
$$r^{d-\alpha} \leq C_2 \alpha \int_{B_r(x) \cap U} \int_{\mathbb{R}^d \setminus U} \frac{\mathrm{d}y\,\mathrm{d}z}{|y-z|^{\alpha+d}}.$$

It thus follows from the first estimate for $\alpha \in [1/2, 1)$ and the second estimate for $\alpha \in (0, 1/2]$ that
$$r^{d-\alpha} \leq C_3 \alpha(1-\alpha) \int_{B_r(x) \cap U} \int_{\mathbb{R}^d \setminus U} \frac{\mathrm{d}y\,\mathrm{d}z}{|y-z|^{\alpha+d}},$$
for every $\alpha \in (0,1)$, where $C_3 = 2\max\{C_1, C_2\}$.

By Wiener's covering lemma, one can extract a countable family of balls $\big(B_{5r(x_i)}(x_i)\big)_{i \in \mathbb{N}}$ with $x_i \in U$ which covers $U$ and is such that the balls $\big(B_{r(x_i)}(x_i)\big)_{i \in \mathbb{N}}$ are disjoint. We thus have that
$$\sum_{i=0}^{\infty} (5r(x_i))^{d-\alpha} \leq C_4 \alpha(1-\alpha) \sum_{i=0}^{\infty} \int_{B_{r(x_i)}(x_i) \cap U} \int_{\mathbb{R}^d \setminus U} \frac{\mathrm{d}y\,\mathrm{d}z}{|y-z|^{\alpha+d}}$$
$$\leq C_4 \alpha(1-\alpha) \int_U \int_{\mathbb{R}^d \setminus U} \frac{\mathrm{d}y\,\mathrm{d}z}{|y-z|^{\alpha+d}} = \frac{C_4}{2} \alpha(1-\alpha) P_\alpha(U).$$

The family $\big(B_{5r(x_i)}(x_i)\big)_{i \in \mathbb{N}}$ thus satisfies the required properties. $\square$

The proof of the coarea formula involving the fractional perimeter $P_\alpha$ is based on the straightforward observation that
$$P_\alpha(A) = [\chi_A]_{W^{\alpha,1}(\mathbb{R}^d)},$$
for every $\alpha \in (0,1)$ and every Borel set $A \subset \mathbb{R}^d$.

**Lemma 4.3.** *For every $u \in W^{\alpha,1}(\mathbb{R}^d)$, we have*
$$[u]_{W^{\alpha,1}(\mathbb{R}^d)} = \int_{-\infty}^{\infty} P_\alpha(\{u > t\})\,\mathrm{d}t.$$

*Proof.* For every $y, z \in \mathbb{R}^d$, we have
$$|u(y) - u(z)| = \int_{-\infty}^{\infty} |\chi_{\{u>t\}}(y) - \chi_{\{u>t\}}(z)| \, \mathrm{d}t.$$

Thus, for $y \neq z$,
$$\frac{|u(y) - u(z)|}{|y - z|^{\alpha+d}} = \int_{-\infty}^{\infty} \frac{|\chi_{\{u>t\}}(y) - \chi_{\{u>t\}}(z)|}{|y - z|^{\alpha+d}} \, \mathrm{d}t.$$

The conclusion follows integrating with respect to $y$ and $z$, and applying Fubini's Theorem. $\square$

*Proof of Theorem 1.3.* Given $t > 0$ and $\varphi \in C_c^\infty(\mathbb{R}^d)$, we apply the fractional Boxing inequality with $U = \{|\varphi| > t\}$. Since $\varphi$ is continuous and has compact support, $U$ is an open bounded subset of $\mathbb{R}^d$, hence there exists a sequence of balls $(B_{r_i}(x_i))_{i \in \mathbb{N}}$ that covers $\{|\varphi| > t\}$ and satisfies
$$\sum_{i=0}^{\infty} r_i^{d-\alpha} \leq C\alpha(1-\alpha) P_\alpha(\{|\varphi| > t\}).$$

It thus follows from the definition of the Hausdorff content that
$$\mathcal{H}_\infty^{d-\alpha}(\{|\varphi| > t\}) \leq \omega_{d-\alpha} C\alpha(1-\alpha) P_\alpha(\{|\varphi| > t\}).$$

Hence,
$$\int_{\mathbb{R}^d} |\varphi| \, \mathrm{d}\mathcal{H}_\infty^{d-\alpha} = \int_0^\infty \mathcal{H}_\infty^{d-\alpha}(\{|\varphi| > t\}) \, \mathrm{d}t$$
$$\leq \omega_{d-\alpha} C\alpha(1-\alpha) \int_0^\infty P_\alpha(\{|\varphi| > t\}) \, \mathrm{d}t.$$

By the fractional coarea formula and the Lipschitz continuity of the absolute-value function,
$$\int_0^\infty P_\alpha(\{|\varphi| > t\}) \, \mathrm{d}t = [|\varphi|]_{W^{\alpha,1}(\mathbb{R}^d)} \leq [\varphi]_{W^{\alpha,1}(\mathbb{R}^d)}.$$

Combining both inequalities, the conclusion follows. $\square$

To deduce the classical Boxing inequality as the limit of the fractional one as $\alpha \to 1$, one first takes a covering of an open set $U \subset \mathbb{R}^d$ of finite perimeter such that
$$\sum_{i=0}^{\infty} r_i^{d-\alpha} \leq C(1-\alpha) P_\alpha(U).$$

Since the sequence of balls does not depend on $\alpha$, we can take $\alpha \to 1$ in the left-hand side. By the second limit in (1.3), the right-hand side converges to $C_1 \operatorname{Per}(U)$ for some constant $C_1 > 0$ independent of $U$, and so one deduces Theorem 1.1. It thus follows from the definition of the Hausdorff content that
$$\mathcal{H}_\infty^{d-1}(U) \leq \sum_{i=0}^{\infty} \omega_{d-1} r_i^{d-1} \leq C_1 \operatorname{Per}(U).$$

Now, given a function $\varphi \in C_c^\infty(\mathbb{R}^d)$, by Sard's lemma the open set $\{|\varphi| > t\}$ is bounded and smooth and so has finite perimeter for almost every $t > 0$. The estimate above applied with $U = \{|\varphi| > t\}$ gives

$$\mathcal{H}_\infty^{d-1}(\{|\varphi| > t\}) \leq C_1 \operatorname{Per}(\{|\varphi| > t\}),$$

for almost every $t > 0$. Integrating both sides with respect to $t$ and using the classical coarea formula, one deduces that

$$\begin{aligned}
(4.3) \quad \int_{\mathbb{R}^d} |\varphi| \, d\mathcal{H}_\infty^{d-1} &\leq C_1 \int_0^\infty \operatorname{Per}(\{|\varphi| > t\}) \, dt \\
&= C_1 \int_{\mathbb{R}^d} |\nabla \varphi| = C_1 [\varphi]_{W^{1,1}(\mathbb{R}^d)}.
\end{aligned}$$

*Proof of Corollary 1.5 and* (1.10), *assuming Theorem 1.7.* As mentioned to in the Introduction, one can argue these results by duality, though we here give a direct argument which relies only upon elementary calculus.

For every open set $U \subset \mathbb{R}^d$, we have

$$\left(\frac{|U|}{\omega_d}\right)^{\frac{d-\alpha}{d}} \leq \frac{\mathcal{H}_\infty^{d-\alpha}(U)}{\omega_{d-\alpha}},$$

which follows from a covering argument and the concavity of the function $s \in [0, \infty) \mapsto s^{\frac{d-\alpha}{d}}$. Applying this estimate with $U = \{|\varphi| > t\}$ and $t > 0$, from Theorems 1.3 and 1.7 we then get

$$\begin{aligned}
(4.4) \quad \int_0^\infty |\{|\varphi| > t\}|^{\frac{d-\alpha}{d}} \, dt &\leq \delta(\alpha, d) \int_0^\infty \mathcal{H}_\infty^{d-\alpha}(\{|\varphi| > t\}) \, dt \\
&\leq \delta(\alpha, d) \, C'(d) \operatorname{dist}(\alpha, \mathbb{N}) [\varphi]_{W^{\alpha,1}(\mathbb{R}^d)},
\end{aligned}$$

where $\delta(\alpha, d) := \omega_d^{\frac{d-\alpha}{d}} / \omega_{d-\alpha}$ is bounded from above independently of $\alpha$ and $d$; see Remark 4.4 below.

Estimate (4.4) yields the continuous embedding of the fractional Sobolev space $W^{\alpha,1}(\mathbb{R}^d)$ into the Lorentz space $L^{\frac{d}{d-\alpha},1}(\mathbb{R}^d)$, which is known to be stronger than the Sobolev embedding in the Lebesgue space $L^{\frac{d}{d-\alpha}}(\mathbb{R}^d)$. We recall such an argument for the sake of completeness (see e.g. [31, Lemma 1.3.5/1]): Since the function $t \mapsto |\{|\varphi| > t\}|$ is nonincreasing, for almost every $t > 0$ one has

$$\frac{d}{dt}\left(\int_0^t |\{|\varphi| > s\}|^{\frac{d-\alpha}{d}} \, ds\right)^{\frac{d}{d-\alpha}} = \frac{d}{d-\alpha}\left(\int_0^t |\{|\varphi| > s\}|^{\frac{d-\alpha}{d}} \, ds\right)^{\frac{d}{d-\alpha}-1} |\{|\varphi| > t\}|^{\frac{d-\alpha}{d}}$$

$$\geq \frac{d}{d-\alpha} t^{\frac{d}{d-\alpha}-1} |\{|\varphi| > t\}|.$$

By the Fundamental Theorem of Calculus for absolutely continuous functions and Cavalieri's principle it follows that

(4.5)
$$\left(\int_0^\infty |\{|\varphi| > t\}|^{\frac{d-\alpha}{d}} \, dt\right)^{\frac{d}{d-\alpha}} \geq \frac{d}{d-\alpha} \int_0^\infty t^{\frac{d}{d-\alpha}-1} |\{|\varphi| > t\}| \, dt = \int_{\mathbb{R}^d} |\varphi|^{\frac{d}{d-\alpha}}.$$

Combining (4.4) and (4.5), we deduce the fractional Sobolev inequality. □

**Remark 4.4.** One shows using basic properties of the Gamma function that

$$\delta(\alpha, d) = \frac{\omega_d^{\frac{d-\alpha}{d}}}{\omega_{d-\alpha}} \leq 1,$$

independently of $\alpha$ and $d$. Indeed, the function $x \mapsto \log(\Gamma(x+1))$ is convex and non-decreasing on the interval $[0, \infty)$ (see e.g. Theorem 1.9 in [9]). Since $\log(\Gamma(1)) = 0$, we thus have that the map

$$x \in (0, \infty) \longmapsto \frac{\log(\Gamma(x+1))}{x}$$

is also non-decreasing. This fact implies that

$$s \in [0, d] \longmapsto \omega_s^{1/s} = \frac{\pi^{1/2}}{(\Gamma(\frac{s}{2}+1))^{1/s}}$$

is a nonincreasing function by computing $\log(\omega_s^{1/s})$. Hence $\omega_d^{\frac{d-\alpha}{d}} \leq \omega_{d-\alpha}$ as claimed.

We now move to the case of a smooth connected bounded open subset $\Omega \subset \mathbb{R}^d$. Given $\alpha \in (0, 1)$ and an open subset $U \subset \Omega$, define the relative perimeter $P_\alpha(U, \Omega)$ by

$$P_\alpha(U, \Omega) := 2 \int_U \int_{\Omega \setminus U} \frac{\mathrm{d}y \, \mathrm{d}z}{|y-z|^{\alpha+d}} = [\chi_U]_{W^{\alpha,1}(\Omega)}.$$

**Proposition 4.5.** *Let $\Omega \subset \mathbb{R}^d$ be a smooth connected bounded open set and $\gamma \in (0, 1)$. There exists a constant $C > 0$ depending on $\gamma$, $d$, and $\Omega$ such that*

$$\mathcal{H}_\infty^{d-\alpha}(U) \leq C(1-\alpha) P_\alpha(U, \Omega),$$

*for every open subset $U \subset \Omega$ with $|U|/|\Omega| \leq \gamma$ and every $\alpha \in (0, 1)$.*

*Proof.* Since $(1+\gamma)/2 \in (|U|/|\Omega|, 1)$, for every $x \in U$ by the Intermediate Value Theorem there exists $r \in (0, \operatorname{diam}\Omega)$ such that

(4.6) $$\frac{|B_r(x) \cap U|}{|B_r(x) \cap \Omega|} = \frac{1+\gamma}{2}.$$

Observe that such a choice depends on $x$ and $U$ but not on $\alpha$.

*Claim.* There exists $\delta = \delta(\Omega) > 0$ such that, for every $\alpha \in (0, 1)$, every open subset $U \subset \Omega$ with $|U|/|\Omega| \leq \gamma$ and every $x \in U$, if $r > 0$ satisfies (4.6) and $r \leq \delta$, then

$$r^{d-\alpha} \leq C'(1-\alpha) \int_{B_r(x) \cap U} \int_{(B_r(x) \cap \Omega) \setminus U} \frac{\mathrm{d}y \, \mathrm{d}z}{|y-z|^{\alpha+d}},$$

for some constant $C' = C'(\Omega, d) > 0$.

*Proof of the Claim.* Observe that for every $\alpha \in (0,1)$ and $U \subset \Omega$ one has

$$\int_{B_r(x) \cap U} \int_{(B_r(x) \cap \Omega) \setminus U} \frac{\mathrm{d}y\, \mathrm{d}z}{|y-z|^{\alpha+d}} \geq \frac{1}{(2r)^{\alpha+d}} |B_r(x) \cap U| |(B_r(x) \cap \Omega) \setminus U|.$$

By the choice of the radius $r$,

$$|B_r(x) \cap U| = \frac{1+\gamma}{2} |B_r(x) \cap \Omega| \quad \text{and} \quad |(B_r(x) \cap \Omega) \setminus U| = \frac{1-\gamma}{2} |B_r(x) \cap \Omega|.$$

By smoothness of $\Omega$, there exists $c > 0$ such that

(4.7) $$|B_r(x) \cap \Omega| \geq c\, |B_r(x)|,$$

for every $x \in \Omega$ and $r \leq \operatorname{diam}\Omega$. We thus have

$$\int_{B_r(x) \cap U} \int_{(B_r(x) \cap \Omega) \setminus U} \frac{\mathrm{d}y\, \mathrm{d}z}{|y-z|^{\alpha+d}} \geq \frac{c_1}{r^{\alpha+d}} \frac{1-\gamma^2}{4} c^2 |B_r(x)|^2 = c_2\, r^{d-\alpha},$$

where $c_2 > 0$ depends on $\gamma$ and $\Omega$.

Now, if the claim is false, then there exist sequences

(a) $\delta_n \to 0$,
(b) $C'_n \to +\infty$,
(c) $(\alpha_n)_{n \in \mathbb{N}}$ with $\alpha_n \in (0,1)$,
(d) $(U_n)_{n \in \mathbb{N}}$ with $U_n \subset \Omega$ open and $|U_n|/|\Omega| \leq \gamma$,
(e) $x_n \to x$ with $x_n \in U_n$ and $x \in \overline{\Omega}$,
(f) $r_n \to 0$ with $r_n \leq \delta_n$ and $|B_{r_n}(x_n) \cap U_n|/|B_{r_n}(x_n) \cap \Omega| = (1+\gamma)/2$,

for which one has

$$r_n^{d-\alpha_n} > C'_n (1-\alpha_n) \int_{B_{r_n}(x_n) \cap U_n} \int_{(B_{r_n}(x_n) \cap \Omega) \setminus U_n} \frac{\mathrm{d}y\, \mathrm{d}z}{|y-z|^{\alpha_n+d}}.$$

Combined with the preceding integral inequality, we conclude firstly that $\alpha_n \to 1$, and secondly

$$\lim_{n \to \infty} (1-\alpha_n) \frac{P_{\alpha_n}\big(B_{r_n}(x_n) \cap U_n, B_{r_n}(x_n) \cap \Omega\big)}{r_n^{d-\alpha_n}} = 0.$$

First we observe that the limit point $x$ must belong to $\partial\Omega$, since otherwise for sufficiently large $n$, $B_{r_n}(x_n) \subset \Omega$, which would lead to a contradiction in light of Lemma 4.1.

Now consider the rescaled open sets

$$A_n := B_1 \cap \frac{U_n - x_n}{r_n} \quad \text{and} \quad V_n := B_1 \cap \frac{\Omega - x_n}{r_n}.$$

They satisfy in particular

(4.8) $$\frac{|A_n|}{|V_n|} = \frac{1+\gamma}{2}, \qquad \frac{|V_n|}{|B_1|} \geq c,$$

and

(4.9) $$\lim_{n \to \infty} (1-\alpha_n) P_{\alpha_n}(A_n, V_n) = 0.$$

Passing to a subsequence if necessary, the sequence of compact sets $(\overline{V_n})_{n\in\mathbb{N}}$ converges with respect to the Hausdorff distance to a compact set $K \subset \overline{B_1}$. In view of the smoothness of $\Omega$, $K$ is the closed ball $\overline{B_1}$ itself or the closed ball intersected with an affine half-space. In particular, the interior $D := \text{int}\, K$ is connected and we also have

(4.10) $$\lim_{n\to\infty} |V_n| = |D|.$$

Given a smooth connected open subset $\omega \Subset D$, there exist $\epsilon > 0$ and $m \in \mathbb{N}$ such that $\omega + B_\epsilon \subset V_n$, for every $n \geq m$. We can further take $\epsilon > 0$ sufficiently small so that $\omega + B_\epsilon$ is also smooth. Since

$$[\chi_{A_n}]_{W^{\alpha_n,1}(\omega+B_\epsilon)} \leq [\chi_{A_n}]_{W^{\alpha_n,1}(V_n)} = P(A_n, V_n),$$

it follows from (4.9) that

(4.11) $$(1-\alpha_n)[\chi_{A_n}]_{W^{\alpha_n,1}(\omega+B_\epsilon)} \leq C_1 \quad \text{for } n \geq m.$$

By a compactness result of Bourgain, Brezis and Mironescu [11, Theorem 4], since $\omega + B_\epsilon$ is smooth, a sequence of functions $(\chi_{A_n})_{n\in\mathbb{N}}$ that satisfies (4.11) is relatively compact in $L^1(\omega + B_\epsilon)$. Passing to a subsequence if necessary, we may thus assume that $(\chi_{A_n})_{n\in\mathbb{N}}$ converges in $L^1(\omega + B_\epsilon)$ and almost everywhere in $\omega + B_\epsilon$ to a function $f$; in particular $f$ must be a characteristic function.

We now use a trick that has been suggested by E. Stein in the context of a new characterization of Sobolev spaces by Bourgain, Brezis and Mironescu [11, 15]: Given a nonnegative function $\rho \in C_c^\infty(B_\epsilon)$ such that $\int_{\mathbb{R}^d} \rho = 1$, by Jensen's inequality we have

$$(1-\alpha_n)[\rho * \chi_{A_n}]_{W^{\alpha_n,1}(\omega)} \leq (1-\alpha_n)[\chi_{A_n}]_{W^{\alpha_n,1}(\omega+B_\epsilon)} \leq (1-\alpha_n)P_{\alpha_n}(A_n, V_n).$$

Since the sequence $(\rho * \chi_{A_n})_{n\in\mathbb{N}}$ is equi-smooth and (4.9) holds, as $n \to \infty$ we then get

$$\|\nabla(\rho * f)\|_{L^1(\omega)} = 0.$$

Hence, by connectedness of $\omega$ we have that $\rho * f$ is constant in $\omega$, and this fact for every $\rho$ implies that $f$ itself is constant in $\omega$.

We claim that for $\eta > 0$ sufficiently small, if $\omega \Subset D$ is such that

$$|D| - |\omega| \leq \eta,$$

then $f = 1$ almost everywhere in $\omega$ and a contradiction then follows. Indeed, since $A_n \subset V_n$ and $\omega \subset V_n$,

$$|A_n \cap \omega| = |A_n| - |A_n \setminus \omega| \geq |A_n| - (|V_n| - |\omega|).$$

By (4.8), we have $|A_n| \geq \lambda$ for every $n \in \mathbb{N}$, with $\lambda := \frac{1+\gamma}{2}c\,|B_1|$. Take $\widetilde{m} \in \mathbb{N}$ such that $|D| \geq |V_n| + \lambda/2$ for $n \geq \widetilde{m}$. Thus,

$$|A_n \cap \omega| \geq \frac{\lambda}{2} - (|D| - |\omega|) \geq \frac{\lambda}{2} - \eta \quad \text{for every } n \geq \widetilde{m}.$$

In view of the pointwise convergence of $(\chi_{A_n})_{n \in \mathbb{N}}$ in $\omega$ we deduce that for $\eta < \lambda/2$, we have $f = 1$ on a subset of positive measure, whence almost everywhere in $\omega$. By (4.10) and the $L^1$ convergence of $(\chi_{A_n})_{n \in \mathbb{N}}$,

$$\frac{|\omega|}{|D|} = \lim_{n \to \infty} \frac{|A_n \cap \omega|}{|V_n|}.$$

Therefore,

$$1 - \frac{\eta}{|D|} \le \frac{|\omega|}{|D|} = \lim_{n \to \infty} \frac{|A_n \cap \omega|}{|V_n|} \le \lim_{n \to \infty} \frac{|A_n|}{|V_n|} = \frac{1 + \gamma}{2}.$$

We thus have a contradiction if we start with $\eta < \min\{\frac{\lambda}{2}, \frac{1-\gamma}{2}|D|\}$. This concludes the proof of the inequality. $\square$

Using Wiener's covering lemma, we cover $U$ with countably many balls $(B_{5r(x_i)}(x_i))_{i \in \mathbb{N}}$ with $x_i \in U$ and such that the balls $(B_{r(x_i)}(x_i))_{i \in \mathbb{N}}$ are disjoint. Since $U$ is fixed, we do not make explicit the dependence of $r$ on $U$. Write $\mathbb{N} = I_1 \cup I_2$ as a disjoint union, where $I_1$ denotes the set of indices $i \in \mathbb{N}$ such that $r(x_i) \le \delta$, where $\delta = \delta(\Omega)$ is given by the Claim. Proceeding as in the proof of Theorem 1.2, we have

$$(4.12) \qquad \sum_{i \in I_1} (r(x_i))^{d-\alpha} \le C_1 (1 - \alpha) \int_U \int_{\Omega \setminus U} \frac{\mathrm{d}y \, \mathrm{d}z}{|y - z|^{\alpha+d}}.$$

We now estimate the sum over $I_2$. By (4.7),

$$|B_{r(x)}(x) \cap U| = \frac{1 + \gamma}{2} |B_{r(x)}(x) \cap \Omega| \ge c'(r(x))^d.$$

Since $r(x_i) > \delta$ for every $i \in I_2$ and the sets $B_{r(x_i)}(x_i) \cap U$ are disjoint, in this case we have

$$\sum_{i \in I_2} (r(x_i))^{d-\alpha} \le \frac{1}{\delta^\alpha} \sum_{i \in I_2} (r(x_i))^d \le \frac{1}{c'\delta^\alpha} \sum_{i \in I_2} |B_{r(x_i)}(x_i) \cap U| \le \frac{1}{c'\delta^\alpha} |U|.$$

Since $\Omega$ is smooth and connected, we can apply the counterpart of the fractional Poincaré inequality (4.1) on the connected open set $\Omega$ (see Corollary 2.5 and Theorem 1.3 in [39]) to get

$$|U||\Omega \setminus U| \le C_2(1 - \alpha) \int_U \int_{\Omega \setminus U} \frac{\mathrm{d}y \, \mathrm{d}z}{|y - z|^{\alpha+d}}.$$

Since $|\Omega \setminus U|/|\Omega| \ge 1 - \gamma$ and $\delta$ depends only on $\Omega$, we then get

$$(4.13) \qquad \sum_{i \in I_2} (r(x_i))^{d-\alpha} \le C_3(1 - \alpha) \int_U \int_{\Omega \setminus U} \frac{\mathrm{d}y \, \mathrm{d}z}{|y - z|^{\alpha+d}}.$$

Combining (4.12) and (4.13), we deduce the estimate for $\mathcal{H}^{d-\alpha}_\infty(U)$. The proof of Proposition 4.5 is complete. $\square$

From Proposition 4.5 we deduce the Poincaré-Wirtinger trace inequality on domains:

*Proof of Theorem 1.6.* Let $\varphi \in C^\infty(\overline{\Omega})$ be such that $\int_\Omega \varphi = 0$. We recall that

(4.14) $$\|\varphi\|_{L^1(\Omega,\,\mathrm{d}\mu)} \leq 2\|\varphi - a\|_{L^1(\Omega,\,\mathrm{d}\mu)},$$

for every $a \in \mathbb{R}$. This is well-known and follows from the observation that

$$|a|\,\mu(\Omega) = \left|\int_\Omega (\varphi - a)\,\mathrm{d}\mu\right| \leq \|\varphi - a\|_{L^1(\Omega,\,\mathrm{d}\mu)},$$

and then apply the triangle inequality to deduce (4.14). In view of (4.14), it thus suffices to prove that

$$\|\varphi - a\|_{L^1(\Omega,\,\mathrm{d}\mu)} \leq C_1(1-\alpha)\|\varphi\|_{W^{\alpha,1}(\Omega)},$$

for some $a \in \mathbb{R}$. For this purpose, we consider two cases:

(i) there exists $a \in \mathbb{R}$ such that $|\{\varphi = a\}| \geq \frac{1}{3}|\Omega|$,
(ii) $|\{\varphi = a\}| < \frac{1}{3}|\Omega|$, for every $a \in \mathbb{R}$.

Assuming that (i) holds, then for every $t > 0$ we have $|\{|\varphi - a| > t\}| \leq \frac{2}{3}|\Omega|$. From the assumption $\mu \leq \mathcal{H}_\infty^{d-\alpha}$ and the relative Boxing inequality with $U = \{|\varphi - a| > t\}$ we have

$$\mu(\{|\varphi - a| > t\}) \leq \mathcal{H}_\infty^{d-\alpha}(\{|\varphi - a| > t\}) \leq C_2\,P_\alpha(\{|\varphi - a| > t\}, \Omega).$$

Integrating both sides with respect to $t$ over $(0, \infty)$ and using the fractional coarea formula with respect to $\Omega$ we get

$$\int_0^\infty \mu(\{|\varphi - a| > t\})\,\mathrm{d}t \leq C_2[|\varphi - a|]_{W^{\alpha,1}(\Omega)} \leq C_2[\varphi]_{W^{\alpha,1}(\Omega)}.$$

We have the conclusion using Cavalieri's principle in the left-hand side.

We now assume that (ii) is satisfied. In this case, the function

$$t \in \mathbb{R} \longmapsto \frac{|\{\varphi < t\}| - |\{\varphi > t\}|}{|\Omega|}$$

has jump discontinuities by less that $1/3$. In addition, it equals $-1$ as $t \to -\infty$ and $1$ as $t \to +\infty$. Hence, it achieves some value in the interval $[0, 1/3]$: There exists $a \in \mathbb{R}$ such that

$$0 \leq \frac{|\{\varphi < a\}| - |\{\varphi > a\}|}{|\Omega|} \leq \frac{1}{3}.$$

Since

$$\frac{|\{\varphi < a\}| + |\{\varphi > a\}|}{|\Omega|} \leq 1,$$

we thus have that $|\{\varphi < a\}| \leq \frac{2}{3}|\Omega|$ and $|\{\varphi > a\}| \leq \frac{2}{3}|\Omega|$. It then follows from the assumption on $\mu$ and the relative Boxing inequality that

$$\mu(\{\varphi - a < -t\}) \leq \mathcal{H}_\infty^{d-\alpha}(\{\varphi - a < -t\}) \leq C_2\,P_\alpha(\{\varphi - a < -t\}, \Omega),$$

for every $t > 0$. Similarly,

$$\mu(\{\varphi - a > t\}) \leq C_2\,P_\alpha(\{\varphi - a > t\}, \Omega),$$

from which we deduce that

$$\mu(\{|\varphi - a| > t\}) \leq C_2\Big(P_\alpha(\{(\varphi - a)^- > t\}, \Omega) + P_\alpha(\{(\varphi - a)^+ > t\}, \Omega)\Big)$$

Integrate both sides with respect to $t$ over $(0, \infty)$. Using Cavalieri's principle and the relative coarea formula we get

$$\|\varphi - a\|_{L^1(\Omega, \mathrm{d}\mu)} \leq C_2\Big([(\varphi - a)^-]_{W^{\alpha,1}(\Omega)} + [(\varphi - a)^+]_{W^{\alpha,1}(\Omega)}\Big) \leq 2C_2[\varphi]_{W^{\alpha,1}(\Omega)},$$

and the conclusion follows. $\square$

The proof of Theorem 1.7 relies on the following estimate by D. Adams involving the Riesz potential $I_\alpha$, which is defined as

$$I_\alpha(z) = \frac{\gamma(\alpha, d)}{|z|^{d-\alpha}},$$

where

(4.15) $$\gamma(\alpha, d) := \frac{\Gamma(\frac{d-\alpha}{2})}{\pi^{d/2}\, 2^\alpha\, \Gamma(\frac{\alpha}{2})}.$$

We sketch the proof of the lemma below to keep track of the dependence of the constant.

**Lemma 4.6.** *Let $k \in \{1, \ldots, d-1\}$ and $s \in (0, d-k)$. There exists a constant $C = C(d) > 0$ such that*

$$\int_{\mathbb{R}^d} I_k * f \, \mathrm{d}\mathcal{H}^s_\infty \leq \frac{C}{d - (s+k)} \int_{\mathbb{R}^d} f \, \mathrm{d}\mathcal{H}^{s+k}_\infty,$$

*for every nonnegative function $f \in C_c^0(\mathbb{R}^d)$.*

*Proof.* By Adams' estimate of the Riesz potential [1, Theorem 6], if $s \in (0, d-1)$, then for every nonnegative locally finite Borel measure $\mu$ in $\mathbb{R}^d$ such that $\mu \leq \mathcal{H}^s_\infty$ and for every Borel set $A \subset \mathbb{R}^d$ one has

$$\int_A I_1 * \mu \leq C_1 \mathcal{H}^{s+1}_\infty(A);$$

see also [40] where the computation on the bottom of p. 294 shows that $C_1 \leq C'/(d - (s+1))$. By Fubini's theorem and Cavalieri's principle, one thus gets

(4.16)
$$\int_{\mathbb{R}^d} I_1 * f \, \mathrm{d}\mu = \int_{\mathbb{R}^d} f\, I_1 * \mu = \int_0^\infty \left(\int_{\{f>t\}} I_1 * \mu\right) \mathrm{d}t$$
$$\leq C_1 \int_0^\infty \mathcal{H}^{s+1}_\infty(\{f > t\}) \, \mathrm{d}t = C_1 \int_{\mathbb{R}^d} f \, \mathrm{d}\mathcal{H}^{s+1}_\infty.$$

To conclude the case $k = 1$, one takes the supremum of the integral in the left-hand side with respect to all nonnegative measures $\mu \leq \mathcal{H}^s_\infty$. For this purpose, one relies on the sublinearity of the Choquet integral for strongly subadditive capacities. The Choquet integral with respect to the dyadic

Hausdorff content $\widehat{\mathcal{H}}_\infty^s$ defined by (3.9) is strongly subadditive, hence by Choquet's theorem the map

$$\varphi \in C_c^0(\mathbb{R}^d) \longmapsto \int_{\mathbb{R}^d} \varphi^+ \, d\widehat{\mathcal{H}}_\infty^s,$$

is sublinear [8, pp. 247–248; 17, Section 54.2]. It thus follows from the Hahn-Banach theorem that the supremum of the left-hand side of (4.16) over non-negative measures $\mu \leq \widehat{\mathcal{H}}_\infty^s$ is comparable to

$$\int_{\mathbb{R}^d} I_1 * f \, d\widehat{\mathcal{H}}_\infty^s,$$

which yields the conclusion of the lemma when $k = 1$ since $\widehat{\mathcal{H}}_\infty^s \sim \mathcal{H}_\infty^s$.

We now proceed by induction on $k$ using the semi-group property of the Riesz potential, which gives $I_k = I_{k-1} * I_1$. To this end, assume the conclusion holds with the Riesz potential $I_{k-1}$, where $k \geq 2$, and we use the fact that $I_k * f = I_{k-1} * (I_1 * f)$. Although $I_1 * f$ need not belong to $C_c^0(\mathbb{R}^d)$, by nonnegativity of $f$ and lower semicontinuity of $I_1 * f$, one finds a non-decreasing sequence $(g_n)_{n \in \mathbb{N}}$ in $C_c^0(\mathbb{R}^d)$ that converges pointwise to $I_1 * f$. Then,

$$\int_{\mathbb{R}^d} I_{k-1} * g_n \, d\mathcal{H}_\infty^s \leq \frac{C_2}{d - (s + k - 1)} \int_{\mathbb{R}^d} g_n \, d\mathcal{H}_\infty^{s+k-1}$$
$$\leq \frac{C_2}{d - (s + k - 1)} \int_{\mathbb{R}^d} I_1 * f \, d\mathcal{H}_\infty^{s+k-1},$$

One can now justify the limit in this inequality as $n \to \infty$ with the Increasing Set Lemma (3.10) for the dyadic Hausdorff content. We then get by semi-group property of the Riesz potential,

$$\int_{\mathbb{R}^d} I_k * f \, d\mathcal{H}_\infty^s = \int_{\mathbb{R}^d} I_{k-1} * (I_1 * f) \, d\mathcal{H}_\infty^s \leq \frac{C_3}{d - (s + k - 1)} \int_{\mathbb{R}^d} I_1 * f \, d\mathcal{H}_\infty^{s+k-1},$$

for every $s \in (0, d - k)$. Applying now the estimate for the Riesz potential $I_1$, we obtain

$$\int_{\mathbb{R}^d} I_k * f \, d\mathcal{H}_\infty^s \leq \frac{C_3}{d - (s + k - 1)} \cdot \frac{C_4}{d - (s + k)} \int_{\mathbb{R}^d} f \, d\mathcal{H}_\infty^{s+k}$$
$$\leq C_3 \cdot \frac{C_4}{d - (s + k)} \int_{\mathbb{R}^d} f \, d\mathcal{H}_\infty^{s+k},$$

which is the inequality we wanted to prove. $\square$

We recall that a Borel measure $\mu$ in $\mathbb{R}^d$ belongs to the Morrey space $M^p(\mathbb{R}^d)$ with $p \geq 1$ if

$$\|\mu\|_{M^p(\mathbb{R}^d)} := \sup_{B_r(x) \subset \mathbb{R}^d} \frac{|\mu|(B_r(x))}{r^{d(p-1)/p}} < \infty.$$

Observe for example that $L^p(\mathbb{R}^d)$ is included in $M^p(\mathbb{R}^d)$ by Hölder's inequality in the sense that the measure $\mu = g \, dx$ belongs to $M^p(\mathbb{R}^d)$ for every $g \in L^p(\mathbb{R}^d)$.

The following lemma from [37] involves the Morrey space $M^d(\mathbb{R}^d)$ and is used in the proof of Theorem 1.7. We sketch the proof to emphasize its connection with the classical Boxing inequality:

**Lemma 4.7.** *For every $\nu \in M^d(\mathbb{R}^d)$, there exists $Y \in L^\infty(\mathbb{R}^d, \mathbb{R}^d)$ such that*
$$-\operatorname{div} Y = \nu \quad \text{in the sense of distributions in } \mathbb{R}^d,$$
*with $\|Y\|_{L^\infty(\mathbb{R}^d)} \leq C\|\nu\|_{M^d(\mathbb{R}^d)}$.*

*Proof.* By linearity of the equation, we may focus on the case where the measure $\nu$ is nonnegative. The assumption $\nu \in M^d(\mathbb{R}^d)$ is then equivalent to $\nu \leq C_1 \|\nu\|_{M^d(\mathbb{R}^d)} \mathcal{H}^{d-1}_\infty$ for some constant $C_1 > 0$ depending on $d$. It follows from the functional formulation (4.3) of the classical Boxing inequality that
$$\left| \int_{\mathbb{R}^d} \varphi \, d\nu \right| \leq C_1 \|\nu\|_{M^d(\mathbb{R}^d)} \int_{\mathbb{R}^d} |\varphi| \, d\mathcal{H}^{d-1}_\infty \leq C_2 \|\nu\|_{M^d(\mathbb{R}^d)} \|\nabla \varphi\|_{L^1(\mathbb{R}^d)},$$
for every $\varphi \in C_c^\infty(\mathbb{R}^d)$. By the Hahn-Banach theorem, the map
$$\nabla \varphi \longmapsto \int_{\mathbb{R}^d} \varphi \, d\nu$$
has a continuous extension in $L^1(\mathbb{R}^d, \mathbb{R}^d)$. The conclusion follows from the Riesz representation theorem. $\square$

*Proof of Theorem 1.7.* The case $\alpha \in (0, 1]$ has already been covered by Theorem 1.3 and (4.3), while the case $\alpha = 0$ follows from the fact that $\mathcal{H}^d_\infty$ coincides with Lebesgue's outer measure in $\mathbb{R}^d$. The remaining case is then $\alpha \in (1, d]$. Given $k \in \{1, \ldots, d-1\}$, we take $\alpha \in (k, k+1]$ and argue differently according to whether $\alpha$ is close to $k$ or $k+1$:

*Case 1.* $\alpha \in [k + \frac{1}{2}, k+1]$.

For every $\varphi \in C_c^\infty(\mathbb{R}^d)$, we have $|\varphi| \leq C_1 \, I_k * |\nabla^k \varphi|$. By monotonicity of the Choquet integral and Adams' integral estimate with $s = d - \alpha$, we have
(4.17)
$$\int_{\mathbb{R}^d} |\varphi| \, d\mathcal{H}^{d-\alpha}_\infty \leq C_1 \int_{\mathbb{R}^d} I_k * |\nabla^k \varphi| \, d\mathcal{H}^{d-\alpha}_\infty \leq \frac{C_2}{\alpha - k} \int_{\mathbb{R}^d} |\nabla^k \varphi| \, d\mathcal{H}^{d-\alpha+k}_\infty.$$
Assuming $\alpha - k \in (0, 1)$, the fractional Boxing inequality gives
$$\int_{\mathbb{R}^d} |\nabla^k \varphi| \, d\mathcal{H}^{d-\alpha+k}_\infty \leq C_3(\alpha - k)(k + 1 - \alpha)[\nabla^k \varphi]_{W^{\alpha-k,1}(\mathbb{R}^d)}.$$
Combining both estimates, we get the conclusion since $k + 1 - \alpha = \operatorname{dist}(\alpha, \mathbb{N})$. When $\alpha = k + 1$, the classical Boxing inequality implies (4.3). Thus,
$$\int_{\mathbb{R}^d} |\nabla^k \varphi| \, d\mathcal{H}^{d-1}_\infty \leq C_4 \|\nabla(\nabla^k \varphi)\|_{L^1(\mathbb{R}^d)} = C_4 \|\nabla^{k+1} \varphi\|_{L^1(\mathbb{R}^d)},$$
which inserted in (4.17) gives the inequality with integer order $k + 1$.

*Case 2.* $\alpha \in (k, k + \frac{1}{2})$.

For any nonnegative Borel measure $\mu$ in $\mathbb{R}^d$ such that $\mu \le \mathcal{H}_\infty^{d-\alpha}$, by Adams' estimate of Riesz potentials in Morrey spaces we have

$$I_{\alpha-1} * \mu \le C_1 \, \mathcal{H}_\infty^{d-1},$$

by identification of the function $I_{\alpha-1} * \mu$ with the measure $I_{\alpha-1} * \mu \, \mathrm{d}x$. In terms of Morrey spaces, $I_{\alpha-1} * \mu \in M^d(\mathbb{R}^d)$. We next observe that for every $\varphi \in C_c^\infty(\mathbb{R}^d)$ we have $\varphi = I_{\alpha-k} * [(-\Delta)^{\frac{\alpha-k}{2}} \varphi]$ and

$$\left|(-\Delta)^{\frac{\alpha-k}{2}} \varphi\right| \le C_2 \, I_{k-1} * \left|\nabla^{k-1}[(-\Delta)^{\frac{\alpha-k}{2}} \varphi]\right| = C_2 \, I_{k-1} * \left|(-\Delta)^{\frac{\alpha-k}{2}} \nabla^{k-1} \varphi\right|,$$

for some constant $C_2 > 0$ depending on $d$. By the semi-group property of the Riesz potential we thus have the pointwise estimate

$$|\varphi| \le C_2 \, I_{\alpha-1} * \left|(-\Delta)^{\frac{\alpha-k}{2}} \nabla^{k-1} \varphi\right|.$$

Integrating this estimate with respect to $\mu$ and applying Fubini's theorem we get

$$\int_{\mathbb{R}^d} |\varphi| \, \mathrm{d}\mu \le C_2 \int_{\mathbb{R}^d} \left|(-\Delta)^{\frac{\alpha-k}{2}} \nabla^{k-1} \varphi\right| I_{\alpha-1} * \mu.$$

Applying Lemma 4.7 with $\nu = I_{\alpha-1} * \mu \, \mathrm{d}x$ we find $Y \in L^\infty(\mathbb{R}^d, \mathbb{R}^d)$ such that

$$-\operatorname{div} Y = I_{\alpha-1} * \mu \quad \text{in } \mathbb{R}^d$$

and $\|Y\|_{L^\infty(\mathbb{R}^d)} \le C_3$, with a constant depending only on $d$. Therefore,

$$\int_{\mathbb{R}^d} |\varphi| \, \mathrm{d}\mu \le C_2 \int_{\mathbb{R}^d} \nabla \left|(-\Delta)^{\frac{\alpha-k}{2}} \nabla^{k-1} \varphi\right| \cdot Y \le C_4 \int_{\mathbb{R}^d} \left|(-\Delta)^{\frac{\alpha-k}{2}} \nabla^k \varphi\right|.$$

Since this estimate holds for every $\mu \le \mathcal{H}_\infty^{d-\alpha}$, taking the supremum of the left-hand side with respect to $\mu$ and applying Choquet's theorem as in the proof of Lemma 4.6, we get

$$\int_{\mathbb{R}^d} |\varphi| \, \mathrm{d}\mathcal{H}_\infty^{d-\alpha} \le C_5 \int_{\mathbb{R}^d} \left|(-\Delta)^{\frac{\alpha-k}{2}} \nabla^k \varphi\right|.$$

Next, by the integral representation of the operator $(-\Delta)^{\frac{\alpha-k}{2}}$, we have

$$\int_{\mathbb{R}^d} \left|(-\Delta)^{\frac{\alpha-k}{2}} \nabla^k \varphi\right| = c_{d,\alpha-k} \int_{\mathbb{R}^d} \left|\int_{\mathbb{R}^d} \frac{\nabla^k \varphi(y) - \nabla^k \varphi(z)}{|y-z|^{\alpha-k+d}} \, \mathrm{d}y\right| \mathrm{d}z$$

$$\le c_{d,\alpha-k} \left[\nabla^k \varphi\right]_{W^{\alpha-k,1}(\mathbb{R}^d)},$$

where

$$c_{d,\alpha-k} := (\alpha - k) \frac{2^{\alpha-k-1} \Gamma(\frac{d+\alpha-k}{2})}{\pi^{\frac{d}{2}} \Gamma(1 - \frac{\alpha-k}{2})}.$$

The result is thus demonstrated, since $c_{d,\alpha-k} = O(\alpha - k)$ as $\alpha \to k^+$. $\square$

We now establish the equivalence between the homogeneous $(\alpha, 1)$-capacity and the $\mathcal{H}_\infty^{d-\alpha}$ Hausdorff content.

*Proof of Theorem 2.1.* Given a set $A \subset \mathbb{R}^d$ and a covering $(B_{r_i}(x_i))_{i \in \mathbb{N}}$ of $A$, by countable subadditivity of the capacity we have

$$\operatorname{Cap}_{\alpha,1}(A) \leq \sum_{i=0}^{\infty} \operatorname{Cap}_{\alpha,1}(B_{r_i}(x_i)).$$

Hence,

$$\operatorname{Cap}_{\alpha,1}(A) \leq \frac{\operatorname{Cap}_{\alpha,1}(B_1)}{\omega_{d-\alpha}} \sum_{i=0}^{\infty} \omega_{d-\alpha} r_i^{d-\alpha},$$

and the estimate $\operatorname{Cap}_{\alpha,1} \leq C_1 \mathcal{H}_\infty^{d-\alpha}$ then follows by taking the infimum in the right-hand side.

To get the reverse estimate, we begin with a compact subset $K \subset \mathbb{R}^d$. Observe that for every nonnegative function $\varphi \in C_c^\infty(\mathbb{R}^d)$ such that $\varphi > 1$ in $K$, by monotonicity of the Hausdorff content and the Chebyshev inequality we have

$$\mathcal{H}_\infty^{d-\alpha}(K) \leq \mathcal{H}_\infty^{d-\alpha}(\{\varphi > 1\}) \leq \int_0^\infty \mathcal{H}_\infty^{d-\alpha}(\{\varphi > t\}) \, dt.$$

From Theorem 1.7 we then get

$$\mathcal{H}_\infty^{d-\alpha}(K) \leq C_2 [\varphi]_{W^{\alpha,1}(\mathbb{R}^d)}.$$

The conclusion in this case follows by taking the infimum with respect to $\varphi$.

The remaining of the proof follows a usual argument based on the inner and outer regularity of set functions: We next consider an open set $U \subset \mathbb{R}^d$. For every compact subset $K \subset U$, by the previous case applied to $K$ and the definition of the capacity on open sets we have

$$\mathcal{H}_\infty^{d-\alpha}(K) \leq C_2 \operatorname{Cap}_{\alpha,1}(K) \leq C_2 \operatorname{Cap}_{\alpha,1}(U).$$

We now take a non-decreasing sequence of compact subsets $(K_n)_{n \in \mathbb{N}}$ such that $\bigcup_{n=0}^{\infty} K_n = U$. Since the dyadic Hausdorff content $\widehat{\mathcal{H}}_\infty^{d-\alpha}$ defined by (3.9) is strongly subadditive, by the Increasing Set Lemma we have

$$\widehat{\mathcal{H}}_\infty^{d-\alpha}(U) = \lim_{n \to \infty} \widehat{\mathcal{H}}_\infty^{d-\alpha}(K_n).$$

Using the previous estimate with $K = K_n$ and the equivalence $\mathcal{H}_\infty^{d-\alpha} \sim \widehat{\mathcal{H}}_\infty^{d-\alpha}$, as $n \to \infty$ we deduce that

$$\mathcal{H}_\infty^{d-\alpha}(U) \leq C_3 \operatorname{Cap}_{\alpha,1}(U).$$

Finally, given any $A \subset \mathbb{R}^d$, by monotonicity of the Hausdorff content and the previous inequality we have

$$\mathcal{H}_\infty^{d-\alpha}(A) \leq \mathcal{H}_\infty^{d-\alpha}(U) \leq C_3 \operatorname{Cap}_{\alpha,1}(U),$$

for every open superset $U \supset A$. The conclusion follows from the definition of $\operatorname{Cap}_{\alpha,1}(A)$ as the infimum of the capacity in the right-hand side. $\square$

## 5. The Fractional Capacity of D. Adams

Setting aside the issue of the continuity of the estimate at the integers, we now show how our Theorems 1.7 and 2.1 can be deduced from the ideas in Adams' paper [2]; see also [3]. The following theorem has been stated by Adams in the case of integers, though implicit in his work is the estimate for $\alpha \in (0, d)$:

**Proposition 5.1.** *Let $\alpha \in (0, d)$. There exists a constant $A = A(\alpha, d) > 0$ such that*
$$\int_{\mathbb{R}^d} |f| \, d\mathcal{H}^{d-\alpha}_\infty \leq A \|(-\Delta)^{\frac{\alpha-k}{2}} |\nabla^k f|\|_{H^1(\mathbb{R}^d)}.$$
*for every $f \in \mathcal{S}_{00}$, where $k = \lfloor \alpha \rfloor$.*

Here, $\mathcal{S}_{00}$ is the space of Schwartz functions in $\mathbb{R}^d$ all of whose moments are zero, and $H^1(\mathbb{R}^d)$ is the real Hardy space. The constant $A$ given by the proof of Proposition 5.1 satisfies
$$A \sim \|I_\alpha\|_{\mathcal{L}(M^{\frac{d}{\alpha}}, \mathrm{BMO})},$$
where the right-hand side denotes the norm of the continuous linear mapping $\mu \mapsto I_\alpha * \mu$ from the Morrey space of measures $M^{\frac{d}{\alpha}}(\mathbb{R}^d)$ into the space of functions of bounded mean oscillation $\mathrm{BMO}(\mathbb{R}^d)$; see [1].

*Proof of Proposition 5.1.* We begin with the pointwise estimate
$$|f| \leq C_1 I_k * |\nabla^k f|.$$
Applying Fubini's theorem, for any nonnegative Borel measure $\mu$ in $\mathbb{R}^d$ we find
$$\int_{\mathbb{R}^d} |f| \, d\mu \leq C_1 \int_{\mathbb{R}^d} I_k * |\nabla^k f| \, d\mu = C_1 \int_{\mathbb{R}^d} [(-\Delta)^{\frac{\alpha-k}{2}} |\nabla^k f|] I_\alpha * \mu.$$
Then by $H^1$-BMO duality (see [22]),
$$\int_{\mathbb{R}^d} |f| \, d\mu \leq C_1 \|(-\Delta)^{\frac{\alpha-k}{2}} |\nabla^k f|\|_{H^1(\mathbb{R}^d)} \|I_\alpha * \mu\|_{\mathrm{BMO}(\mathbb{R}^d)}$$
$$\leq C_1 \|(-\Delta)^{\frac{\alpha-k}{2}} |\nabla^k f|\|_{H^1(\mathbb{R}^d)} \|I_\alpha\|_{\mathcal{L}(M^{\frac{d}{\alpha}}, \mathrm{BMO})} \|\mu\|_{M^{\frac{d}{\alpha}}(\mathbb{R}^d)}.$$
Taking the supremum in the left-hand side over all nonnegative measures $\mu \in M^{\frac{d}{\alpha}}(\mathbb{R}^d)$ such that $\|\mu\|_{M^{\frac{d}{\alpha}}(\mathbb{R}^d)} \leq 1$ and applying Choquet's theorem as in the proof of Lemma 4.6 above we deduce that
$$\int_{\mathbb{R}^d} |f| \, d\mathcal{H}^{d-\alpha}_\infty \leq C_2 \|(-\Delta)^{\frac{\alpha-k}{2}} |\nabla^k f|\|_{H^1(\mathbb{R}^d)} \|I_\alpha\|_{\mathcal{L}(M^{\alpha/d}, \mathrm{BMO})},$$
which gives the conclusion. □

**Remark 5.2.** A slight variation in the proof which relies on Strichartz characterization of these Hardy-Sobolev spaces (see [42]) shows one has the inequality

$$\int_{\mathbb{R}^d} |f| \, d\mathcal{H}^{d-\alpha}_\infty \leq A' \, \|(-\Delta)^{\frac{\alpha}{2}} f\|_{H^1(\mathbb{R}^d)}. \tag{5.1}$$

for any $\alpha \in (0, d)$.

The fact that there is some constant $C = C(\alpha, d) > 0$ for which the conclusion of Theorem 1.7 holds is then a consequence of Proposition 5.1 and the following inequality applied to $g = |\nabla^k f|$:

$$\|(-\Delta)^{\frac{\alpha-k}{2}} g\|_{H^1(\mathbb{R}^d)} \leq B \, [g]_{W^{\alpha-k,1}(\mathbb{R}^d)}, \tag{5.2}$$

where $\alpha - k \in (0, 1)$ and $B = B(\alpha, d) > 0$ is a constant. The validity of (5.2) follows from the observation that

$$\|(-\Delta)^{\frac{\alpha-k}{2}} g\|_{H^1(\mathbb{R}^d)} = \tilde{c}_{d,\alpha-k} \int_{\mathbb{R}^d} \left| \int_{\mathbb{R}^d} \frac{g(y) - g(z)}{|y-z|^{\alpha-k+d}} \cdot \frac{y-z}{|y-z|} \, dy \right| dz$$

$$+ c_{d,\alpha-k} \int_{\mathbb{R}^d} \left| \int_{\mathbb{R}^d} \frac{g(y) - g(z)}{|y-z|^{\alpha-k+d}} \, dy \right| dz \leq B \, [g]_{W^{\alpha-k,1}(\mathbb{R}^d)}.$$

In particular, one finds the following version of our estimate:

$$\int_{\mathbb{R}^d} |\varphi| \, d\mathcal{H}^{d-\alpha}_\infty \leq AB \, [\varphi]_{W^{\alpha,1}(\mathbb{R}^d)}, \tag{5.3}$$

for every $\varphi \in C^\infty_c(\mathbb{R}^d)$ and $\alpha \in (0, d)$. However, we claim that the constant $AB$ obtained this fashion cannot tend to 0 as $\operatorname{dist}(\alpha, \mathbb{N})$ at the integers. More precisely, we now show that the constant $A = A(\alpha, d)$ obtained in Proposition 5.1 and any constant $B = B(\alpha, d)$ for which inequality (5.2) holds one has

$$\lim_{\alpha \to n} \frac{A(\alpha, d) B(\alpha, d)}{|\alpha - n|} = +\infty$$

for $n \in \{1, 2, \ldots, d\}$. To see this, let us first establish

**Proposition 5.3.** *The operator norm of the Riesz potential $I_\alpha : M^{\frac{d}{\alpha}}(\mathbb{R}^d) \to \operatorname{BMO}(\mathbb{R}^d)$ satisfies*

$$\|I_\alpha\|_{\mathcal{L}(M^{\frac{d}{\alpha}}, \operatorname{BMO})} \geq C\alpha,$$

*for every $\alpha \in (0, d)$ and some constant $C = C(d) > 0$.*

*Proof.* Observe that for every ball $B_r(x) \subset \mathbb{R}^d$,

$$\int_{B_r(x)} \frac{1}{|z|^\alpha} \, dz \leq \int_{B_r(0)} \frac{1}{|z|^\alpha} \, dz = \frac{\operatorname{Per}(B_1)}{d - \alpha} r^{d-\alpha}.$$

Thus, the measure $\mu_\alpha := I_{d-\alpha} \, dz$ belongs to the Morrey space $M^{\frac{d}{\alpha}}(\mathbb{R}^d)$ and

$$\|\mu_\alpha\|_{M^{\frac{d}{\alpha}}(\mathbb{R}^d)} = \gamma(d - \alpha, d) \frac{\operatorname{Per}(B_1)}{d - \alpha}.$$

One also finds
$$I_\alpha * \mu_\alpha(x) = \frac{\mathrm{Per}\,(B_1)}{(2\pi)^d} \log \frac{1}{|x|} = \frac{1}{2^{d-1}\pi^{\frac{d}{2}}\Gamma(\frac{d}{2})} \log \frac{1}{|x|};$$

see e.g. p. 50 in [27]. Here one should pay attention that the normalization in [27] is such that the Fourier transform of his Riesz potentials have Fourier symbol $|\xi|^{-\alpha}$, while the definition of our constant $\gamma(\alpha, d)$ ensures the symbol is $(2\pi|\xi|)^{-\alpha}$, and hence the additional factor of $(2\pi)^d$ in the denominator.

In particular, the quantity $\|I_\alpha * \mu_\alpha\|_{\mathrm{BMO}(\mathbb{R}^d)}$ is independent of $\alpha$. Hence,

$$\|I_\alpha\|_{\mathcal{L}(M^{\frac{d}{\alpha}},\mathrm{BMO})} \geq \frac{\|I_\alpha * \mu_\alpha\|_{\mathrm{BMO}(\mathbb{R}^d)}}{\|\mu_\alpha\|_{M^{\frac{d}{\alpha}}(\mathbb{R}^d)}} = c(d)\frac{d-\alpha}{\gamma(d-\alpha,d)}.$$

By the explicit formula of $\gamma(d-\alpha,d)$ given by (4.15) and standard properties of the Gamma function, we have $(d-\alpha)/\gamma(d-\alpha,d) = O(\alpha)$ and we validate our claim. $\square$

Next, we show

**Proposition 5.4.** *For every $n \in \{0, 1, \ldots, d\}$ one has*
$$\lim_{\alpha \to n} \frac{B(\alpha, d)}{|\alpha - n|} = +\infty.$$

*Proof.* Given $k \in \{0, 1, \ldots, d-1\}$ and $\alpha \in (k, k+1)$, suppose that $B(\alpha, d) > 0$ satisfies
$$\|(-\Delta)^{\frac{\alpha-k}{2}}\varphi\|_{H^1(\mathbb{R}^d)} \leq B(\alpha, d)\,[\varphi]_{W^{\alpha-k,1}(\mathbb{R}^d)}$$
for every $\varphi \in C_c^\infty(\mathbb{R}^d)$ (or more generally for Lipschitz functions with compact support). Then inserting the factor of $(\alpha - k)(k + 1 - \alpha)$ we have

$$\|(-\Delta)^{\frac{\alpha-k}{2}}\varphi\|_{H^1(\mathbb{R}^d)} \leq \frac{B(\alpha, d)}{(\alpha-k)(k+1-\alpha)}(\alpha-k)(k+1-\alpha)[\varphi]_{W^{\alpha-k,1}(\mathbb{R}^d)}.$$

Now as $\alpha \to (k+1)^-$, one obtains
$$\|(-\Delta)^{\frac{1}{2}}\varphi\|_{H^1(\mathbb{R}^d)} \leq \left[\liminf_{\alpha \to (k+1)^-} \frac{B(\alpha,d)}{k+1-\alpha}\right]\|\nabla\varphi\|_{L^1(\mathbb{R}^d)},$$

and so it is not possible that the limit on the right-hand side stays bounded, since this would yield a false embedding. Similarly, as $\alpha \to k^+$ one finds
$$\|\varphi\|_{H^1(\mathbb{R}^d)} \leq \left[\liminf_{\alpha \to k^+} \frac{B(\alpha,d)}{\alpha-k}\right]\|\varphi\|_{L^1(\mathbb{R}^d)},$$

which cannot hold with finite constant. The result is thus demonstrated. $\square$

Let us now relate this discussion to the capacity $\mathrm{Cap}_{\alpha,1}$ introduced in Section 2. Observe that Remark 5.2 implies

(5.4) $$\mathcal{H}_\infty^{d-\alpha}(K) \leq CR'_\alpha(K),$$

for every compact set $K \subset \mathbb{R}^d$, where

$$R'_\alpha(K) := \inf \{ \|(-\Delta)^{\frac{\alpha}{2}} f\|_{H^1(\mathbb{R}^d)} : f \in \mathcal{S}_{00},\ f \geq 1 \text{ on } K \}.$$

Adams introduced in [2] the capacity

$$R_\alpha(K) := \inf \{ \|g\|_{H^1(\mathbb{R}^d)} : g \in \mathcal{S}_{00},\ I_\alpha * g \geq 1 \text{ on } K \}.$$

In general one cannot restrict the computation of this capacity to non-negative test functions, and as a result $R_\alpha$ may fail to be countably subadditive. However, we observe that

**Proposition 5.5.** *Let $\alpha \in (0, d)$. For every compact set $K \subset \mathbb{R}^d$,*

$$R_\alpha(K) = R'_\alpha(K).$$

*Proof.* First let us remark that the space $\mathcal{S}_{00}$ is closed with respect to the operator $I_\alpha$ and its inverse $(-\Delta)^{\frac{\alpha}{2}}$. Therefore, if $g \in \mathcal{S}_{00}$ with $I_\alpha * g \geq 1$, the function $f = I_\alpha * g \in \mathcal{S}_{00}$ satisfies $f \geq 1$ and restricting oneself to the consideration of such $f$ one finds

$$\inf \{ \|g\|_{H^1(\mathbb{R}^d)} : I_\alpha * g \geq 1 \} \leq \inf \{ \|(-\Delta)^{\frac{\alpha}{2}} f\|_{H^1(\mathbb{R}^d)} : f \geq 1 \},$$

and thus $R_\alpha(K) \leq R'_\alpha(K)$. Conversely, for any $f \in \mathcal{S}_{00}$ with $f \geq 1$, restricting consideration in $R'_\alpha$ to functions $g = (-\Delta)^{\frac{\alpha}{2}} f \in \mathcal{S}_{00}$ one finds the reverse inequality. □

Then the lower bound (5.4), the embedding implied by inequality (5.2), and the straightforward estimate (2.3) yield the chain of inequalities

$$\mathcal{H}^{d-\alpha}_\infty(K) \leq C R_\alpha(K) \leq C' \operatorname{Cap}_{\alpha,1}(K) \leq C'' \mathcal{H}^{d-\alpha}_\infty(K),$$

so that these quantities are all equivalent. In this way we obtain a different approach to Theorem 2.1, and also the following result of Adams:

**Proposition 5.6.** *Let $\alpha \in (0, d)$. For every compact set $K \subset \mathbb{R}^d$,*

$$R_\alpha(K) \sim \mathcal{H}^{d-\alpha}_\infty(K).$$

## ACKNOWLEDGEMENTS

Part of this work was written while the first author (ACP) was visiting NCTU with support from the Taiwan Ministry of Science and Technology through the Mathematics Research Promotion Center. He warmly thanks the Department of Applied Mathematics for the invitation and hospitality. The second author (DS) is supported by the Taiwan Ministry of Science and Technology under research grant 105-2115-M-009-004-MY2.

AUGUSTO C. PONCE
UNIVERSITÉ CATHOLIQUE DE LOUVAIN
INSTITUT DE RECHERCHE EN MATHÉMATIQUE ET PHYSIQUE
CHEMIN DU CYCLOTRON 2, L7.01.02
1348 LOUVAIN-LA-NEUVE, BELGIUM
*E-mail address*: `Augusto.Ponce@uclouvain.be`

DANIEL SPECTOR
NATIONAL CHIAO TUNG UNIVERSITY
DEPARTMENT OF APPLIED MATHEMATICS
HSINCHU, TAIWAN

NATIONAL CENTER FOR THEORETICAL SCIENCES
NATIONAL TAIWAN UNIVERSITY
NO. 1 SEC. 4 ROOSEVELT RD.
TAIPEI, 106, TAIWAN
*E-mail address*: `dspector@math.nctu.edu.tw`